\documentclass[11pt]{amsart}
\usepackage{amsfonts,amssymb,amsmath}
\usepackage[latin1]{inputenc}
\usepackage[dvips]{epsfig}
\usepackage{graphicx}
\newtheorem{theorem}{Theorem}[section]
\newtheorem{lemma}[theorem]{Lemma}
\newtheorem{assumption}[theorem]{Assumption}
\newtheorem{corollary}[theorem]{Corollary}
\newtheorem{proposition}[theorem]{Proposition}
\newtheorem{definition}[theorem]{Definition}
\newtheorem{remark}[theorem]{Remark}
\newcommand{\bd}{\begin{displaymath}}
\newcommand{\ed}{\end{displaymath}}
\def\beglem{\begin{lemma}\sl }
\def\endlem{\end{lemma}}
\def\begthm{\begin{theorem}\sl }
\def\endthm{\end{theorem}}
\def\begprop{\begin{proposition} \sl}
\def\endprop{\end{proposition}}
\def\begcor{\begin{corollary}\sl }
\def\endcor{\end{corollary}}
\def\begdef{\begin{definition}\sl}
\def\enddef{\end{definition}}
\def\begproof{ \noindent {\em Proof: \ }}
\def\endproof{\null\hfill {$\Box$}\bigskip}  
\def\beghyp{\begin{assumption} \sl}
\def\endhyp{\end{assumption}}
\def\begrem{\begin{remark}\rm}
\def\endrem{\null\hfill {$\Box$}\end{remark}}

 \newcommand{\RR}{\mathbb R}
 \newcommand{\NN}{\mathbb N} 
 \newcommand{\ZZ}{\mathbb Z}
 \newcommand{\FF}{\mathcal F}

\def\signrd{
\vbox{{\sc Roland Duclous (corresponding author)\par\vspace{3mm}
Université Bordeaux I, \par
UMR CELIA CEA, CNRS et Institut de Mathématiques de Bordeaux \par
351, Cours de la Lib\'{e}ration\par    
F-33405 Talence cedex,
FRANCE\par\vspace{3mm}
e-mail:} duclous@math.u-bordeaux1.fr }
\par\vspace{3mm}
}

\def\signff{
\vbox{{\sc Francis Filbet\par\vspace{3mm}
Universit\'e Lyon,  \par 
Universit\'e Lyon1, CNRS, \par  
UMR 5208 - Institut Camille Jordan, \par 
43, Boulevard du 11 Novembre 1918,\par 
F-69622 Villeurbanne cedex,
FRANCE\par\vspace{3mm}
e-mail:} filbet@math.univ-lyon1.fr }
\par\vspace{3mm}
}

\def\signbd{
\vbox{{\sc Bruno Dubroca\par\vspace{3mm}
Université Bordeaux I, \par
UMR CELIA CEA, CNRS et Institut de Mathématiques de Bordeaux \par
351, Cours de la Lib\'{e}ration\par    
F-33405 Talence cedex,
FRANCE\par\vspace{3mm}
e-mail:} dubroca@math.u-bordeaux1.fr }
\par\vspace{3mm}
}

\def\signvt{
\vbox{{\sc Vladimir Tikhonchuk\par\vspace{3mm}
Université Bordeaux I, \par
UMR CELIA CEA, CNRS \par
351, Cours de la Lib\'{e}ration\par    
F-33405 Talence cedex,
FRANCE\par\vspace{3mm}
e-mail:} tikhonchuk@celia.u-bordeaux1.fr }
\par\vspace{3mm}
}

\begin{document}

\signrd
\signff
\signbd
\signvt

\newpage

\title[High order resolution of the Maxwell-Fokker-Planck-Landau model
intended for ICF applications] {High order resolution of the
Maxwell-Fokker-Planck-Landau model intended for ICF applications}

\author{Roland Duclous, Bruno Dubroca, Francis Filbet and Vladimir Tikhonchuk}

\hyphenation{bounda-ry rea-so-na-ble be-ha-vior pro-per-ties
cha-rac-te-ris-tic}

\begin{abstract}
A high order, deterministic direct numerical method is proposed for
the nonrelativistic $2D_{\bf x} \times 3D_{\bf v}$ Vlasov-Maxwell
system, coupled with Fokker-Planck-Landau type operators. Such a system
is devoted to the modelling of electronic transport and energy
deposition in the general frame of Inertial Confinement Fusion
applications. It describes the kinetics of plasma physics in the
nonlocal thermodynamic equilibrium regime. Strong numerical
constraints lead us to develop specific methods and approaches for
validation, that might be used in other fields where couplings between
equations, multiscale physics, and high dimensionality are
involved. Parallelisation (MPI communication standard) and fast
algorithms such as the multigrid method are employed, that make this
direct approach be computationally affordable for simulations of
hundreds of picoseconds, when dealing with configurations that present
five dimensions in phase space.

\end{abstract}

\maketitle


\medskip
\noindent {\sc AMS subject classifications.}

\tableofcontents

\section{Introduction}

\noindent
In the context of the interaction of intense, short laser pulses with
solid targets \cite{santos}, Inertial Confinement Fusion (ICF) schemes
\cite{bell,tabak}, the energy transport is an important issue. In this
latter field of applications (ICF), it determines the efficiency of
plasma heating and the possibility to achieve the fusion conditions.
The appropriate scales under consideration here are about one hundred
of micrometers for the typical spatial sizes, and one hundred of
picoseconds for the time scales.

Several key features should be accounted for. First of all, in typical
ICF configurations, a significant amount of energetic electrons have a
large mean free path, exceeding the characteristic gradient length of
the temperature and the density: the particles motion exhibits
nonlocal features.

A wide range of collisional regimes should be dealt with to describe
the propagation and the deposit of energetic electrons from the
underdense corona of the target to its dense and compressed core.

The collisions are important even if the beam particles themselves are
collisionless \cite{bell} : these particles, when propagating in a
plasma, trigger a return current that neutralizes the incident
current.  This return current is determined by collisions of thermal,
background electrons.  The structure of the generated electron
distribution function is then often anisotropic, with a strongly
intercorrelated two population structure.  For nonrelativistic laser
intensities, smaller than $10^{18}\,\mbox{W}\,\mbox{cm}^{-2}$, a small
angle description for collisions between the two populations is
well-suited, leading to the classical Fokker-Planck-Landau collision
model.  The Coulomb potential involves a large amount of collisions
with small energy exchanges between particles, so that the Landau form
of the Fokker-Planck operator is required here. Such a configuration
with two counterstreaming beams typically leads to the developement of
microscopic instabilities that can modify strongly the beam
propagation.  We refer to the two-stream and filamentation
instabilities, where the wave vector of the perturbation is
respectively parallel and perpendicular to the incident beam
\cite{bret,bret2}. A self-consistent description of electromagnetic
fields is then required to describe the plasma behaviour and associated
instabilities.  Furthermore in the process of plasma heating, strong
magnetic fields are generated at intensity that can reach a MegaGauss
scale and may affect the energy transport
\cite{braginskii,haines,Kingham}. The sources of magnetic field
generation include on the one hand the effects of the rotational part
of the electronic pressure which is a cross gradient $\nabla T \times
\nabla n$ effect, and on the other hand the exponential growth of
perturbations of anisotropic distribution functions (Weibel
instability). Some electromagnetic processes can be strongly coupled
with nonlocal effects.

The plasma model studied in this paper is based on the nonrelativistic
Vlasov-Maxwell equations, coupled with Fokker-Planck-Landau collision
operators. It gathers the listed requirements at laser intensities
which are relevant for ICF. At higher laser intensities, a
relativistic treatment should be considered \cite{bell,yokota}, and
collision operators with large energy exchanges are required if
secondary fast electron production proves to be non-negligible,
paticularly with dense plasmas.

 There are several numerical methods that treat the Vlasov-Maxwell
model together with Fokker-Planck-Landau type operators. Among them,
the Particle-In-Cell (PIC) methods provide satisfying results only in
a limited range of collisional regimes.  Moreover, they suffer from
the ``finite grid instability'', that leads to numerical heating.
Also the statistical noise and the low resolution of the electron
distribution function by PIC solvers lead generally to an inaccurate
treatment of collisions, particularly when dealing with low
temperature and high density plasmas. Another approach consists in the
expansion of the distribution function in Legendre polynomials,
retaining the lowest order terms. However, with this approach, a
strong anisotropy of the distribution function cannot be treated
\cite{bell,kingham0,shoucri}.  Both of these methods are well-suited
in particular regimes but fail at modelling more complicated
situations where a collisionless anisotropic fast electron population
is coupled to a collision dominated thermal population.  To overcome
these difficulties, a spherical harmonic expansion has been proposed,
that proves to be efficient \cite{bell}.  Here we propose a different
approach which consists in approximating the full model by a direct
deterministic numerical method.  It discretizes directly the initial
set of equations and enables to preserve, at a discrete level, the
physical invariants of the model (positivity of the distribution
function, total mass and energy, entropy decreasing behaviour, etc).
Many deterministic schemes of this type have already been considered
for homogeneous Fokker-Planck type operators
\cite{degond,buetcordier,buetcordier1,Pareshi}. The nonhomogeneous
case, that includes the transport part (see \cite{filbet3} for a
comparison between Eulerian Vlasov solvers), involves a large
computational complexity that can only be reduced with fast
algorithms. Multipole expansion \cite{lemou} and multigrid
\cite{buetcordier1} techniques, as well as fast spectral methods
\cite{Pareshi,filbet2}, have been applied to the Landau equation.  For
computational complexity constraints, very few results on the accuracy
of these methods are known in the nonhomogeneous case
\cite{crouseilles,filbet2}, particularly when the coupling with
magnetic fields is considered.

Our starting point for the transport part discretization is a second
order finite volume scheme introduced in \cite{crouseilles}.  Its main
feature is that it preserves exactly the discrete energy, if slope
limiters are not active. We intoduce additional dissipation on these
limiters in order to successfully address the two-stream instability
test case. We will underline the important role of the limitation
procedure for the accuracy, on the second order scheme. This scheme is
compared in this test case with a fourth order MUSCL scheme
\cite{yamamoto}, with a limitation ensuring the positivity of the
distribution function \cite{berthon1}.  A similar approach, with the
introduction of a fourth order scheme for transport to avoid numerical
heating, has already been proposed in the context of PIC solvers
\cite{sentoku}. The discretization of the Maxwell equations is
performed with a Crank-Nicholson method, allowing to have time steps
of the order of the collision time.  It is designed to preserve the
discrete total electromagnetic energy, which is a very important
numerical constraint when considering the coupling of Vlasov and
Maxwell equations for applications aiming at capturing an accurate
energy deposition.  We use for the Landau operator a fast multigrid
technique that proves to be accurate in a wide range of collisional
regimes.  Moreover, the use of domain decomposition techniques and
distributed memory MPI standard on the space domain leads to
affordable computational cost, allowing to treat time dependent
$2D_{\bf x} \times 3D_{\bf v}$ problems. As for the Lorentz
electron-ion collision operator, we insist on discrete symmetry
properties that are important when coupling to the Maxwell equations.

Finally, we propose to validate the numerical method on several
physical test cases.

The paper is organized as follows. First, we present the model and its
properties, then we discuss the numerical schemes for the transport
part, their properties, and propose several numerical tests.  Then the
discretization for the collision operators is treated and we finally
present physical test cases $1D_{x}\times 3D_{\bf v}$ and $2D_{\bf
x}\times 3D_{\bf v}$ that show the accuracy of the present algorithm.

\section{Kinetic model}
Two particle species are considered: ions which are supposed to be
fixed (assuming an electron-ion mass ratio $m_e/m_i \ll 1$), and
electrons for which the evolution is described by a distribution
function $f_e(t,{\bf x},{\bf v})$ where for the more general case
$({\bf x},{\bf v})\,\in\, \Omega \times \RR^3$, with $\Omega \subset
\RR^3$.  The nonrelativistic Vlasov equation with Fokker-Planck-Landau
(FPL) collision operator is given by
\begin{equation} 
  \label{MFP1}
\frac{\partial f_e}{\partial t} \,+\, \nabla_{{\bf x}}\cdot \left(
{\bf v}\, f_e \right) \,+\, \frac{q_e}{m_e}\,\nabla_{{\bf v}}
\cdot\left( ({\bf E} \,+\, {\bf v} \times{\bf B} )\, f_e \right) \,=\,
C_{e,e}(f_e,f_e)\,+\,C_{e,i}(f_e),
\end{equation}
where $q_e=-e$ is the charge of an electron and $m_e$ is the mass of
an electron. On the one hand, electromagnetic fields $({\bf E}, {\bf
B})$ are given by the classical Maxwell system
\begin{equation} 
\label{MFP2}
\left\{
\begin{array}{l}
 \displaystyle{\frac{\partial {\bf E}}{\partial t}\,-\,c^2
\,\nabla_{\bf x} \times {\bf B} \,=\, -\frac{\bf J}{\epsilon_0}}, \\ \, \\
\displaystyle{\frac{\partial {\bf B}}{\partial t}+ \nabla_{\bf x} \times
{\bf E} \,=\, 0,}
\end{array}\right.
\end{equation}
where $\epsilon_0$ represents the permittivity of vacuum, $c$ is the
speed of light. The electric current is given by
$$ {\bf J}(t,{\bf x}) \,=\, q_e \int_{\RR^3} f_e(t,{\bf x},{\bf v})
\,{\bf v\, dv}.
$$ 
Moreover, Maxwell system's is supplemented by Gauss law's
\begin{equation} 
\label{MFP5}
\nabla_{\bf x}\cdot {\bf E} \,=\, \frac{\rho}{\epsilon_0},  \quad
\nabla_{\bf x}\cdot {\bf B} \,=\, 0,
\end{equation}
where $\rho$ is the charge density:
$$ \rho \,=\, q_e\,(n_e-{n_o}) \,=\, q_e\,\left(\int_{\RR^3}
f_e(t,{\bf x},{\bf v})\, {\bf dv} \,-\, {n_0} \right),
$$ and ${n_0}/{Z}$ is the initial ion density.

On the other hand in (\ref{MFP1}), the right hand side represents
collisions between particles, which only act on the velocity variable,
so we drop the ${\bf x}$ variable.  The operator $C_{e,e}(f_e,f_e)$
stands for the electron-electron collision operator whereas
$C_{e,i}(f_e)$ is the electron-ion collision operator
\begin{equation}
 \label{Cee}
  C_{e,e}(f_e,f_e) \,=\, \frac{e^4 \,\ln \Lambda }{8
 \,\pi\,\epsilon_0^2\, m_e^2} \,\nabla_{\bf v} \cdot \left(\int_{\RR^3}
 \Phi({\bf v}-{\bf v}') \left [ f_e({\bf v}') \nabla_{\bf v} f_e({\bf
 v}) - f_e( {\bf v}) \nabla_{\bf v'} f_e({\bf v}') \right ] {\bf d
 v'}\right),
\end{equation}
 whereas $C_{e,i}(f_e)$ is the electron-ion collision operator
\begin{equation}
 \label{Cei}
  C_{e,i}(f_e) \,=\, \frac{Z \,n_0 \,e^4 \ln\Lambda}{8
  \,\pi\,\epsilon_0^2 \,m_e^2} \,\nabla_{\bf v} \cdot \left [
  \Phi({\bf v}) \nabla_{\bf v} f_e ({\bf v}) \right ].
\end{equation}
where $ \ln\Lambda$ is the Coulomb logarithm, which is supposed to
be constant over the domain and $\Phi({\bf u})$ is an operator acting
on the relative velocity ${\bf u}$
\begin{eqnarray} 
  \label{op}
\Phi({\bf u})= \frac{ \|{\bf u}\|^2\, {\rm Id} \,-\, {\bf u}
\otimes {\bf u}}{\|{\bf u}\|^3}.
\end{eqnarray}
As we assume ions to be fixed, the FPL operator can then be simplified
for electron-ion collisions \cite{crouseilles}, and reduced to the Lorentz
approximation. We refer to \cite{balescu1} for a physical derivation.

In this model, the Vlasov equation stands for the invariance of the
distribution function along the particles trajectories affected by the
electric and magnetic fields ${\bf E}$ and ${\bf B}$.  The Vlasov
equation representing the left-hand side in (\ref{MFP1}) is written in
a conservative form, but it can also be written in an equivalent
non-conservative form, while Maxwell equations
(\ref{MFP2})-(\ref{MFP5}) provide with a complete self-consistent
description of electromagnetic fields. The coupling between both is
performed via the Lorentz force term ${\bf E} \,+\, {\bf v} \times
{\bf B} $ in the Vlasov equation, and the current source terms in
Maxwell equations. Furthermore, the FPL operator is used to describe
elastic, binary collisions between charged particles, with the
long-range Coulomb interaction potential.  Classical but important
properties of the system (\ref{MFP1})-(\ref{MFP5}) together with
operators (\ref{Cee}) and (\ref{Cei}), are briefly recalled. For
detailed proofs, we refer to \cite{crouseilles,decoster}.
%We present expressions (\ref{Cee}) and (\ref{Cei}) in the velocity space domain (Restriction to the homogeneous case can be considered here, since space variable ${\bf x}$ only appears as a parameter in the collision operators ) and in the nonrelativistic limit \cite{decoster}.

\subsection{Transport equation under electromagnetic fields}

Let us neglect in this section the collision operators. The
Vlasov-Maxwell system (\ref{MFP1})-(\ref{MFP5}) with a zero right-hand
side is strictly equivalent to (\ref{MFP1})-(\ref{MFP2}) provided
Gauss's laws (\ref{MFP5}) are initially satisfied. This gives a
compatibility condition at initial time.\\ The mass and momentum are
preserved with respect to time for the Vlasov-Maxwell system, {\it
i.e.} system (\ref{MFP1})-(\ref{MFP2}) without collision operators
$$ 
\frac{d}{dt} \int_{\RR^3 \times \RR^3} f_e(t,{\bf x},{\bf v})
\left( \begin{array}{l} 1 \\ {\bf v} \end{array}\right) {\bf dx\,dv}
\,=\, 0, \,\,t \geq 0.
$$ Moreover, conservation of energy can be proved for the
Vlasov-Maxwell system by multiplying equation (\ref{MFP1}) by
$m_e\,\|{\bf v}\|^2/2$ and integrating it in the velocity space. It
gives after an integration by parts
$$ \frac{1}{2}\,\frac{d}{dt} \int_{\RR^3} \left\{\epsilon_0\,\|{\bf
E}(t,{\bf x})\|^2 \,+\, \frac{1}{\mu_0}\|{\bf B}(t,{\bf x})\|^2 \,+\,
\left[ \int_{\RR^3} m_e \|{\bf v}\|^2 f_e(t,{\bf x},{\bf v}){\bf dv}
\right]\, \right\}\,{\bf dx} \,=\,0, \,\, t \geq 0,
$$ with $c^2\,\epsilon_0 \,\mu_0 \,=\, 1$. The Vlasov-Maxwell system
also conserves the kinetic entropy
$$ \frac{d}{dt} H(t) \,=\, \frac{d}{dt}\int_{\RR^3 \times \RR^3}
f_e(t,{\bf x},{\bf v}) \log(f_e(t,{\bf x},{\bf v}))\,{\bf dx\,dv}\,=\,
0, \,\,t \geq 0.
$$

\subsection{Collision operators}

The FPL operator is used to describe binary elastic collisions between
electrons. Its algebraic structure is similar to the Boltzmann
operator, in that it satisfies the conservation of mass, momentum and
energy
$$
\int_{\RR^3} C_{e,e}(f_e,f_e)({\bf v}) \left(
       \begin{array}{c} 1 \\ {\bf v}\\ \|{\bf v}\|^2 
\end{array}\right) {\bf dv}\,=\,0, \,\, t \geq 0.
$$ 
Moreover, the entropy is decreasing with respect to time
$$ 
\frac{dH}{dt}(t) = \frac{d}{dt} \int_{\RR^3} f_e({\bf
v},t)\,\log(f_e({\bf v},t) ) {\bf dv}\, \le\, 0, \,\,t\geq 0.
$$ 
The equilibrium states of the FPL operator, {\it i.e.} the set of
distribution functions in the kernel of $C_{e,e}(f_e,f_e)$, are given
by the Maxwellian distribution functions
$$
{\mathcal M}_{n_e ,{\bf u_e} , T_e} ({\bf v}) \,=\, n_e \left (
\frac{m_e}{2\,\pi\, T_e} \right )^{3/2} \exp \left( -m_e\,\frac{\|{\bf
v}-{\bf u_e} \|^2}{2 \,T_e} \right),
$$
where $n_e $ is the density, ${\bf u_e}$ is the mean velocity and
$T_e$ is the temperature, defined as
$$ 
\left\{
\begin{array}{l}
\displaystyle{
n_e \,=\, \int_{\RR^3} f_e({\bf v}){\bf dv},}
\\
\,
\\
\displaystyle{{\bf u_e} \,=\,\frac{1}{n_e} \int_{\RR^3} f_e({\bf v}){\bf v}{\bf dv},}
\\
\,
\\
\displaystyle{T_e\,=\,\frac{m_e}{3\,n_e}\int_{\RR^3} f_e({\bf v}) \|{\bf v}-{\bf u_e} \|^2 {\bf dv}.}
\end{array}\right.
$$ 
On the other hand, the operator (\ref{Cei}), modelling collisions
between electrons and ions, is a Lorentz operator. It satisfies the
conservation of mass and energy
$$
 \int_{\RR^3} C_{e,i}(f_e)({\bf v})\left(
       \begin{array}{c} 1 \\ \|{\bf v}-{\bf u_e} \|^2 \end{array}\right) \,{\bf dv} \,=\,0. 
$$ 
Moreover, the equilibrium states for this operator are given by the
set of isotropic functions:
$$
{\rm Ker}\left(C_{e,i}\right)\,=\,\left \{\,f_e\in L^1\left((1+\|{\bf v}\|^2){\bf dv}\right)  ,  \quad f_e({\bf v}) \,=\, \phi(z), \,\,\,z\,=\, \|{\bf v} -{\bf u_e}\|^2\, \right \}.
$$ 
Finally, each convex function $\psi$ of $f_e$ is an entropy for $C_{e,i}(f_e)$,
$$ 
\frac{d}{dt} \int_{\RR^3} \psi (f_e ) \,{\bf dv} \,\le\, 0, \,\,t \geq 0.
$$ 

In addition to these properties, we present a symmetry property.  This
property may have some importance, in particular in presence of
magnetic fields. In that case, any break of symmetry due to an
inadequate discretization method could lead to generation of an
artificial magnetic field, via the current source terms, when coupling
with the Maxwell equations.

\begin{proposition}
\label{prop:01}
\noindent
 If $f_e$ has the following symmetry property with respect to the
 direction $k$ at time $t_0$
\begin{equation}
 \label{symmetryprop0}
f_e(t_0,{\bf v})\,\,=\,\,f_e(t_0,{\bf v}^k),
\end{equation}
with components for
 $$ { \bf v}^k_i \,=\, \left \{
\begin{array}{ll} 
+{\bf v}_i &\textrm{ if } i \neq k, \\ -{\bf v}_i & \textrm{ if } i
= k.
\end{array}
\right .
$$ Then, this symmetry property is preserved with respect to time.
\end{proposition}

\section{Numerical scheme for transport}

We present a finite volume approximation for the Vlasov-Maxwell system
(\ref{MFP1})-(\ref{MFP2}) without collision operators. Indeed, it is
crucial to approximate accurately the transport part of the system to
asses the collective behaviour\footnote[1]{By collective effects, we
denote here the self-consistent interaction of electromagnetic fields
and particles. Some collective effects are also considered in the
collision processes, which make two particles interact via the Coulomb
field. The self-consistent electromagnetic field then screens the long
range Coulomb potential and removes the singularity in the
Fokker-Plank-Landau operator.} of the plasma, that occurs typically at
a shorter scale than the collision processes. We introduce a uniform
$1D$ space discretization $(x_{i+1/2})_{i \in I}$, $I\subset \NN$, of
the interval $(0,L_1)$, in the direction denoted by index $1$. The
associated space variable is denoted by $x_1$. We define the control
volumes $ C_{i,{\bf j}}=(x_{i-1/2},x_{i+1/2})\times({\bf v}_{\bf
j-1/2},{\bf v}_{\bf j+1/2})$, the size of a control volume in one
direction in space $\Delta x$ and velocity $\Delta v $.

The velocity variable ${\bf v}={}^{t}(v_1,v_2,v_3)$ is discretized on
the grid ${\bf v}_{\bf j}\,=\, {\bf j}\,\Delta v \,=\,
{}^{t}(v_{j_1},v_{j_2},v_{j_3})$ with ${\bf
j}={}^{t}(j_1,j_2,j_3)\in\ZZ^3$. Moreover we note ${\bf v}_{\bf
j+1/2}\,=\, {}^{t}(j_1+1/2,j_2+1/2,j_3+1/2)\,\Delta v$. Finally, the
time discretization is defined as $t^n\,=\,n\Delta t$, with $n\in
\NN$.

Let $f_{i,{\bf j}}^n $ be an average approximation of the distribution
function on the control volume $ C_{i,{\bf j}}$ at time $t^n$, that is
$$ f_{i,{\bf j}}^n\,\simeq \,\frac{1}{\Delta x\,\Delta
v^3}\int_{C_{i,{\bf j}}} f(t^n,x,{\bf v}) dx\,{\bf dv}.
$$ Moreover since the discretization is presented in a simple $1D_x$
space geometry, the electromagnetic field has the follownig structure:
${\bf E} = {}^{t}(E_1(t,x_1)),E_2(t,x_1),0)$, ${\bf
B}={}^{t}(0,0,B_3(t,x_1))$. Hence ${}^{t}(E_{1,i}^n,E_{2,i}^n)$ is an
approximation of the electric field ${}^{t}(E_1,E_2)$ whereas
$B_{3,i}^n$ represents an approximation of the magnetic field $B_3$ in
the control volume $(x_{i-1/2},x_{i+1/2})$ at time $t^n$.

\subsection{Second order approximation of a one dimensional transport equation}

For the sake of simplicity, we focus on the discretization of a $1D$
transport equation; the extension to higher dimensions is
straightforward on a grid, without requiring time splitting techniques
between transport terms. In this section, the index $1$ is dropped
both on space and velocity directions, for this simple $1D_x$
geometry.

Let us consider the following equation for $t\geq 0$ and $x \in
(0,L)$,
\begin{equation} 
  \label{equationsimple}
 \frac{\partial f}{\partial t} \,+\, v\,\frac{\partial f}{\partial x}
\,=\, 0 \ ,
\end{equation}
 where the velocity $v > 0$ is given. By symmetry it is possible to
 recover the case when $v$ is negative. In the following we skip the
 velocity variable dependency of the distribution function. Using a
 time explicit Euler scheme and integrating the $1D$ Vlasov equation
 on a control volume $(x_{i-1/2},x_{i+1/2})$, it yields
\begin{equation} 
  \label{contvol0}
f_{i}^{n+1} \,=\, f_{i}^{n} - \frac{\Delta t}{\Delta x} \, \left[
\FF_{i+1/2}^n - \FF_{i-1/2}^n \right],
\end{equation}
where $\FF_{i+1/2}^n$ represents an approximation of the flux $v\,f(t^n,x_{i+1/2})$ at the interface $x_{i+1/2}$. 

The next step consists in approximating the fluxes and to reconstruct
the distribution function.  To this aim, we approximate the
distribution function $f(t^n,x)$ by $f_h(x)$ using a second order
accurate approximation of the distribution function on the interval $
[x_{i-1/2},x_{i+1/2})$, with a reconstruction technique by primitive
\cite{crouseilles}
\begin{equation} 
  \label{slopelim+}
 f_h(x) = f_i^n \,+\, \epsilon_i^{+} \,\frac{(x-x_i)}{\Delta x} (
 f_{i+1}^n \,-\, f_{i}^n ).
\end{equation}
We introduce the limiter
\begin{equation} 
  \label{slopedef+}
 \epsilon_i^{+} = \left \{
\begin{array}{ll} 
0 & \textrm{ if } (f_{i+1}^n - f_{i}^n)\,( f_{i}^n - f_{i-1}^n)\,<\,0,
\\ \min\left(1,\displaystyle{\frac{2 \left ( \|f^0\|_{\infty}-f^n_i
\right ) }{f^n_{i}-f^n_{i+1}}}\right) & \textrm{ if } (f_{i+1}^n -
f_{i}^n) <0 \ , \\ \min \left
(1,\displaystyle{\frac{2\,f_{i}^n}{f_{i+1}^n-f_{i}^n}} \right) &
\textrm{ else,}
\end{array} \right .
\end{equation}
and set $\FF_{i+1/2}^n = v\, f_h(x_{i+1/2})$.  This type of limiter
introduces a particular treatment for extrema. At this price only
(dissipation at extrema), we were able to recover correctly the
two-stream instabililty test case, without oscillations destroying the
salient features of the distribution function structure. Another
choice for the limitation consits in choosing the ``Van Leer's one
parameter family of the minmod limiters'' \cite{Kurganov}
\begin{equation} 
  \label{slopedef+1}
 \epsilon_i^{+} = \mbox{minmod} \left (
 b\frac{(f_{i+1}^n-f_{i}^n)}{\Delta x},
 \frac{(f_{i+1}^n-f_{i-1}^n)}{2\Delta x} ,
 b\frac{(f_{i}^n-f_{i-1}^n)}{\Delta x}\right ) , \, \end{equation}
 where
$$ \mbox{minmod}(x,y,z) \equiv \mbox{max} ( 0 , \mbox{min}(x,y,z))+
 \mbox{min} ( 0 , \mbox{max}(x,y,z)) \ , \quad (x,y,z) \in \RR³ ,
 $$and $b$ is a parameter between $1$ and $2$. We will see on the
 two-stream instability test case the importance of the choice for
 limiters.\\ Finally, this reconstruction ensures the conservation of
 the average and maximum principle on $f_h(x)$ \cite{crouseilles}.

\subsection{Fourth order transport scheme}
 We turn now to a higher order approximation (fourth order MUSCL TVD
scheme) \cite{yamamoto}. This scheme has also been considered in
\cite{berthon1}, in the frame of VFRoe schemes for the shallow water
equations, where the authors proposed an additional limitation. Here
we note that an optimized limitation procedure is possible in our
case, breaking the similar treatment for both right and left
increments, and taking advantage of the structure of the flux in the
nonrelativistic Vlasov equation: the force term does not depend of the
advection variable.\\

For this MUSCL scheme, we only provide
here with the algorithm for the implementation of this scheme and
refer to \cite{berthon1}, \cite{yamamoto} for the derivation procedure
of this scheme.\\
The high order flux at the interface $x_{i+1/2}$, at time $t^n$ reads
\begin{eqnarray}
\displaystyle
 \FF_{i+1/2}^n= \FF \left ( f_{i,r}^n , f_{i+1,l}^n \right ) = \left\{
\begin{array}{c}
 v f_{i,r}^n  \quad \ \ \ \mbox{if}  \quad v>0 \ , \\
 v f_{i+1,l}^n \quad \mbox{if}  \quad v<0 \ .
\end{array}
\right .
\nonumber
\end{eqnarray}
This numerical flux involves the reconstructed states: $\displaystyle
f_{i,r}^n = f_{i}^n + (\Delta f)_i^+$ and $\displaystyle f_{i,l}^n =
f_{i}^n + (\Delta f)_i^- ,$ where $ (\Delta f)_i^{\pm} $ are the
reconstruction increments.\\ An intermediate state $ f_{i }^* $,
defined by $ \displaystyle \frac{1}{3} \left ( f_{i,r}^n + f_{i }^* +
f_{i,l}^n \right ) = f_{i}^n $ si introduced. It is shown in
\cite{berthon1} that the introduction of this intermediate state
preserves, provided the CFL condition is formally divided by three,
the positivity of the distribution function. Following \cite{yamamoto}
and \cite{berthon1}, the fourth order MUSCL reconstruction reads

\vspace{.1in}
\begin{center}
\begin{tabular}{|l|} 
\hline
{\bf Algorithm of reconstruction.} 
\\
\,
\\
Compute 
\\
\hspace{1.cm} $\displaystyle{\left ( \Delta f\right )_i^- = -\frac{1}{6} \left ( 2 \Delta^* \bar{f}_{i-1/2} + \Delta^* \tilde{f}_{i+1/2} \right ), }$
\\
\,
\\
\hspace{1.cm} $\displaystyle{\left ( \Delta f\right )_i^+ = \frac{1}{6} \left ( \Delta^* \bar{f}_{i-1/2} + 2 \Delta^* \tilde{f}_{i+1/2} \right ),}$
\\
where
\\
\hspace{1.cm} $\displaystyle{\Delta^* \bar{f}_{i-1/2} =\mbox{minmod}( \Delta^* f_{i-1/2}, 4 \Delta^* f_{i+1/2} ),}$
\\
\,
\\
\hspace{1.cm} $\displaystyle{\Delta^* \tilde{f}_{i+1/2} =\mbox{minmod}( \Delta^*
f_{i+1/2}, 4 \Delta^* f_{i-1/2} ) }$
\\
and
\\
\hspace{1.cm} $\displaystyle{ \Delta^* f_{i+1/2} = \Delta f_{i+1/2} - \frac{1}{6} \Delta^3 \bar{f}_{i+1/2},}$
\\
\,
\\
\hspace{1.cm} $\displaystyle{ \Delta^3 \bar{f}_{i+1/2} = \Delta \bar{f}^{a}_{i-1/2} - 2 \Delta \bar{f}^b_{i+1/2} + \Delta \bar{f}^c_{i+3/2},}$
\\
with
\\
\hspace{1.cm} $\displaystyle{\Delta \bar{f}^a_{i-1/2} = \mbox{minmod}( \Delta f_{i-1/2}, 2 \Delta f_{i+1/2},2 \Delta f_{i+3/2}),}$
\\
\;
\\
\hspace{1.cm} $\displaystyle{ \Delta \bar{f}^b_{i+1/2} = \mbox{minmod}( \Delta f_{i+1/2}, 2 \Delta f_{i+3/2},2 \Delta f_{i-1/2}),}$
\\
\,
\\
\hspace{1.cm} $\displaystyle{ \Delta \bar{f}^c_{i+3/2} = \mbox{minmod}( \Delta f_{i+3/2}, 2 \Delta f_{i-1/2},2 \Delta f_{i+1/2}),}$
\\
\,
\\
with the notation $\Delta f_{i+1/2}=f_{i+1}-f_i$.
\\
\,
\\
\hline
\end{tabular}
\end{center}

\begin{eqnarray}
\label{minmod2}
&& \mbox{ Reminding that the minmod limiter is given by} \nonumber \\
&& \nonumber \\
&& \mbox{minmod}(x,y) = 
\left\{
\begin{array}{ll}
0,  & \textrm{ if } x\,y \leq 0,
\\
\,
\\
x   & \textrm{ if } |x|\,\leq |y|,
\\
\,
\\
y & \textrm{else,}
\end{array}\right. \nonumber \\
&& \mbox{with $(x,y) \in \mathbb{R}^3$. } \nonumber
\end{eqnarray}

\begin{center}
The limitation proposed in \cite{berthon1} is then applied.\\ It allows to
satisfy the positivity of the reconstructed states.
\end{center}

\vspace{.1in}
\begin{center}
\begin{tabular}{|l|} 
\hline
%\vspace{0.2cm}
{\bf Algorithm for the limitation involving the intermediate state.}
\\
\,
\\
Compute $(\Delta f)_i^{{\rm lim},\pm}$ such that
\\
\hspace{1.cm} $\displaystyle{ 
f_i^n + ( \Delta f )_i^{{\rm lim},-} \ge 0,}$ 
\\
\,
\\
\hspace{1.cm} $\displaystyle{   
f_i^n + ( \Delta f )_i^{{\rm lim},+} \ge 0,}$
\\
and 
\\
\hspace{1.cm} $\displaystyle{f^*_i = f_i^n - ( \Delta f )_i^{{\rm lim},-} - ( \Delta f )_i^{{\rm lim},+} \ge 0.}$
\\
\,
\\
This limitation reads:
\\
\,
\\
\hspace{1.cm} $\displaystyle{
\left \{
\begin{array}{l}
(\Delta f )_i^{{\rm lim},-} \,=\,\theta\, \max\left( ( \Delta f )_i^-,-f_i^n\right),
\\
\,
\\
(\Delta f)_i^{{\rm lim},+} \,=\,\theta \,\max\left( ( \Delta f )_i^+,-f_i^n \right),
\end{array}\right. }$
\\
where
\\
\hspace{1.cm} $\displaystyle{\theta \,=\, \left \{
 \begin{array}{ll}
1,  \quad \mbox{if} \quad \max\left( ( \Delta f )_i^- , -f_i^n \right) \,+\, \max\left( ( \Delta f )_i^+ , -f_i^n \right) \le 0 \ ,
 \\  
 \,
\\
\min \left( 1, \frac{f_i^n}{\max \left( ( \Delta f )_i^- , -f_i^n \right) \,+\, \max \left( ( \Delta f )_i^+ , -f_i^n \right)} \right ) \quad
 \mbox{otherwise.} 
\end{array}
\right. }$
\\
\,
\\
\hline
\end{tabular}
\end{center}

\subsection{Application to the Vlasov-Maxwell system.}

We exactly follow the same idea to design a scheme for the full Vlasov
equation in phase space $({\bf x}, {\bf v})\in \Omega\times\RR^3$.  In
addition, a centered formulation for the electromagnetic fields is
chosen:
\begin{equation} 
  \label{EMfield}
{\bf E}^{n+1/2} \,=\; \frac{1}{2} \left({\bf E}^{n+1} \,+\,{\bf E}^n
\right)   \quad \mbox{and} \quad 
 {\bf B}^{n+1/2} \,=\, \frac{1}{2}\left({\bf B}^{n+1} +{\bf B}^n \right).
\end{equation}
The discretization of the Maxwell equations (\ref{MFP2})-(\ref{MFP5})
is performed via an implicit $\theta$-scheme, with $\theta=1/2$, which
corresponds to the Crank-Nicholson scheme and thus preserves the total
discrete energy. This discretization is presented in a simple $1D$
space geometry. The electric field ${\bf E} = {}^{t}(E_1,E_2,0)$ and
the magnetic field ${\bf B}= {}^{t}(0,0,B_3)$ are collocated data on
the discrete grid. These fields are solution of the system
\begin{equation}
\label{theta}
\left\{\begin{array}{l}
\displaystyle{\frac{E^{n+1}_{1,i}-E^{n}_{1,i}}{\Delta t}\,=\,-
\frac{{J}_{1,i}^n}{\epsilon_0}}, 
\\ \, 
\\
\displaystyle{\frac{E^{n+1}_{2,i}-E^{n}_{2,i}}{\Delta t}
\,+\,c^2\,\frac{B^{n+1/2}_{3,{i+1}}-B^{n+1/2}_{3,{i-1}}}{2\Delta
x}\,=\,- \frac{{J}_{2,i}^n}{\epsilon_0}}, 
\\ 
\, 
\\
\displaystyle{\frac{B^{n+1}_{3,i}-B^{n}_{3,i}}{\Delta t} \,+\,
\frac{E^{n+1/2}_{2,{i+1}}-E^{n+1/2}_{2,{i-1}}}{2\Delta x}\,=\,0}.
\end{array}\right.
\end{equation}  
This scheme is well suited for the electrodynamics situations that are
treated here in the test cases.\\ The approximation for the current in
\eqref{theta} ${J}_{1}^n$ and ${J}_{2}^n$ has been chosen such as
\begin{equation}
\label{closurecurrents}
 {J}_{1,i}^n \,=\, \sum_{{\bf j}\in\ZZ^3} \Delta v^3\, v_{j_1} \,f_{{
i},{\bf j}}^n \quad\textrm{ and }\quad {J}_{2,i}^n \,=\, \sum_{{\bf j}\in\ZZ^3}
\Delta v^3\, v_{j_2} \,f_{{ i},{\bf j}}^n.
\end{equation} 
Unfortunately, these expressions do not preserve the total energy when
slopes limiters are active, but we will show that they have the
important feature to reproduce the discrete two-stream dispersion
relation.

First, we remind discrete properties concerning positivity, mass and
energy conservation \cite{crouseilles} of the second order scheme
(\ref{contvol0})-(\ref{slopelim+}) coupled with
(\ref{EMfield})-(\ref{closurecurrents}), considering now the magnetic
component.

\begin{proposition}
 Let the initial datum $(f_{{\bf i},{\bf j}}^0)_{{\bf i},{\bf
 j}\in\ZZ^3}$ be nonnegative and assume the following $CFL$ type
 condition on the time step
\begin{equation}
\label{CFL}
\Delta t \le C \min\left(\Delta x,\Delta v\right),
\end{equation}  
where $C>0$ is related to the maximum norm of the electric and magnetic fields and the upper bound of the velocity domain. 

Then the scheme (\ref{contvol0})-(\ref{slopelim+}) coupled with
(\ref{EMfield})-(\ref{closurecurrents}), when extended to the infinite
$3D_{\bf x} \times 3D_{\bf v}$ geometry, gives a nonnegative
approximation, preserves total mass and preserves total energy when
slopes limiters are not active on the transport in the velocity
directions
$$ 
\frac{1}{2}  \sum_{{\bf i}\in I}\Delta x^3  \left\{ \epsilon_0\,\|{\bf E}^n_{\bf i}\|^2 \,+\, \frac{1}{\mu_0} \|{\bf B}_{{\bf i}}^n\|^2 \,+\,  m_e\,\left[\sum_{{\bf j}\in\ZZ^3} \Delta v^3\, \|{\bf v}_{{\bf j}}\|^2\,f_{{\bf i},{\bf j}}^n\right] \right\} \,=\, C^0 \ , \ n \in \mathbb{N}.
$$
\end{proposition}

In addition to these properties, we justifiy our choice for the
numerical current thanks to a discrete dispersion relation on the
two-stream instability. In the rest of the section, we drop the index
$1$ on the variables $x_1$, $v_1$, $E_1$ and $J_1$, because the
transport is considered $1D_x \times 1D_v$.
\begin{proposition}
 Consider the second order scheme (\ref{contvol0})-(\ref{slopelim+})
 coupled with (\ref{EMfield})-(\ref{closurecurrents}), when slope
 limiters are not active, to approximate the Vlasov-Amp\`{e}re system
\begin{equation}
\label{simple2stream}
\left\{\begin{array}{l}
\displaystyle{\frac{\partial f}{\partial t} \,+\, v\,\frac{\partial f}{\partial x} \,+\, \frac{q_e}{m_e}\,E\, \frac{\partial f}{\partial v} \,=\, 0,} 
\\
\,
\\
\displaystyle{\frac{\partial E}{\partial t} \,=\, -\frac{J}{\epsilon_0}}.
\end{array}\right.
\end{equation}
Then the definition (\ref{closurecurrents}) for the current $J$
defines a discrete dispersion relation that converges toward the
continuous dispersion relation when $\Delta v$, $\Delta x$ and $\Delta
t$ tend to $0$.
\label{prop:3.2}
\end{proposition}

\begproof The two-stream instability configuration can be fully
 analysed with the Vlasov-Amp\`ere system (\ref{simple2stream})
 extracted from equations (\ref{MFP1})-(\ref{MFP5}). The dispersion
 relation for a perturbation $f^{(1)} \propto e^{i(k\,x\,-\,\omega\,
 t)}$ of an initial equilibrium state $f^{(0)}$, with $\|f^{(1)}\| \ll
 \|f^{(0)}\|$, then reads
\begin{equation}
\label{dispersioncont}
1\,+\,\frac{q_e^2}{\epsilon_0\, m_e} \, \int_{\RR}
\frac{v}{\omega(\omega-k\,v)}\frac{\partial f^{(0)}}{\partial v}\,dv
\,=\,0.
\end{equation}
 Here the crucial point is the discretization on the velocity part of
the phase space, so that we perform a semi-discrete analysis. In the
frame of the discretization (\ref{contvol0})-(\ref{slopelim+}) coupled
with (\ref{EMfield})-(\ref{closurecurrents}), we consider the
semi-discrete scheme approximating (\ref{simple2stream})
\begin{equation}
\label{scheme1D1V1}
\left\{
\begin{array}{l}
\displaystyle{\frac{\partial f}{\partial t} \,+\, v\,\frac{\partial
f}{\partial x} \,+\, \frac{q_e}{m_e}\,E\,
\frac{f_{j+1/2}-f_{j+1/2}}{\Delta v} \,=\, 0,} \\ \, \\
\displaystyle{\frac{\partial E}{\partial t} \,=\, -
\frac{q_e}{\epsilon_0}\, \sum_{j\in\ZZ} \Delta v\, \,v_{j}\, f_{j},}
\end{array}\right.
\end{equation}
with
$$ f_{j+1/2} \,=\, \frac{f_{j+1}+f_{j}}{2},
$$ assuming the slope limiter is not active.  Then we performe a
discrete linearization around an equilibrium state
$$ f_{j} = {f}_{j}^{(0)} \,+\,{f}_{j}^{(1)},
$$ where $\|f^{(1)}\| \ll \|f^{(0)}\|$. Using ${f}_{j}^{(1)} \propto
e^{i(k\,x-\omega t)}$ in (\ref{scheme1D1V1}), it yields
\begin{equation}
\label{schemelimit}
\left\{
\begin{array}{l}
\displaystyle{-i\,(\omega \,-\,k\, v_j) \,f_{j}^{(1)}
\,+\,\frac{q_e}{m_e} \,E^{(1)}
\,\frac{f_{j+1/2}^{(0)}-f_{j-1/2}^{(0)}}{\Delta v} \,=\,0 }, \\ \, \\
\displaystyle{-i\,\omega\, E^{(1)} \,=\, - \frac{q_e}{\epsilon_0}
\,\sum_{j\in\ZZ} \Delta v \, v_j \,f_{j}^{(1)}.}
\end{array}\right.
\end{equation}
These equations lead to the discrete dispersion relation
\begin{equation}
\label{discretedispersionfinal}
1\,+\,\frac{q_e^2}{\epsilon_0 \,m_e} \, \sum_{j\in\ZZ}
\frac{v_{j}}{\omega\,(\omega\,-\,k\,v_j)}\, \left [
\frac{f^{(0)}_{j+1/2} -f^{(0)}_{j-1/2}}{\Delta v} \right ] \Delta v
\,=\, 0.
\end{equation}
We recover the continuous dispersion relation (\ref{dispersioncont})
when passing at the limit $\Delta v\rightarrow 0$. Any other choice
for the discrete current in (\ref{schemelimit}) would introduce an
additional error to the $ O(\Delta v^2)$ error in the relation
dispersion (\ref{discretedispersionfinal}). For instance, choosing
$$
J \,=\, \sum_{j\in\ZZ} \Delta v\, v_j\, f_{j+1/2}
$$
would have lead to the analogous of (\ref{discretedispersionfinal}):
\begin{equation}
\label{schemelimitana}
1\,+\,\frac{q_e^2}{\epsilon_0 \,m_e} \sum_{j\in\ZZ}
\frac{(v_{j}-\Delta v)}{\omega\,(\omega\,-\,k\,v_j)} \,\left [ \frac{
f^{(0)}_{j+1/2} \,-\, f^{(0)}_{j-1/2} }{\Delta v} \right ] \Delta v
\,=\,0,
\end{equation}
which is a ``shifted'' dispersion relation, with a $O(\Delta v)$
 accuracy, compared to the $O(\Delta v^2)$ accuracy on relation
 (\ref{discretedispersionfinal}).  \endproof

\section{Validation of the transport schemes}

We first propose a validation stategy in the linear, collisionless
regime, based on the work of Sartori and Coppa \cite{sc}.  They
performed a transient analysis, and obtain exact solutions of the
periodic Vlasov-Poisson system, in the nonrelativistic and
relativistic regime.

Their approach, relying on Green kernels, is recalled in Appendix
\ref{annexeA}, in the nonrelativistic regime.  A generalization of the
2D periodic relativistic Vlasov-Maxwell system, including magnetic
fields, will be presented in a forthcoming paper. Our objective is to
capture kinetic effects in the linear regime, such as the Landau
damping and the two-stream instability. A semi-analytical solutin is
obtained, with a prescribed accuracy. Moreover, this method allows to
explore wavenumber ranges where other approaches relying on dispersion
relations fail. We recall that classical validations of kinetic
solvers dedicated to plasma physics \cite{crouseilles,bourdiec} are based
on the calculation of the growth rates (instability), or decrease
rates (damping) in the linear regime. Let us show the efficiency of
the semi-analytical method on the two-stream instability test case.

\subsection{Scaling with plasma frequency}

Scaling parameters can be introduced to obtain a dimensionless form of
the Vlasov-Maxwell-Fokker-Planck equations.  The plasma frequency
$\omega_{pe}$, the Debye length $\lambda_D$, the thermal velocity of
electrons $v_{th}$, and the cyclotron frequency $\omega_{ce}$ are
defined as follows
\begin{eqnarray} 
  \label{scaling}
  \omega_{pe} = \sqrt{\frac{n_0 e^2}{\epsilon_0 m_e}}, \quad \lambda_D
  = \sqrt{ \frac{\epsilon_0 \kappa_B T_0}{n_0 e^2}}, \quad v_{th} =
  \sqrt{\frac{\kappa_B T_0}{m_e}}, \quad \omega_{ce} = {\frac{e
  B}{m_e}}.
\end{eqnarray}
These parameters enable us to define dimensionless parameters marked
with tilde.
\begin{itemize}
\item Dimensionless time, space and velocity, respectively:
\begin{equation}
\label{scaling0}
\tilde{t}\,=\,\omega_{pe}\,t, \quad \tilde{x}\,=\;\frac{x}{\lambda_D},
\quad \tilde{v}\,=\,\frac{v}{v_{th}}.
 \end{equation}

\item Dimensionless electric field, magnetic field and distribution
function, respectively

\begin{equation}
  \label{scaling1}
\tilde{E} \,=\, \frac{e\,E}{ m_e \,v_{th} \omega_{pe}}, \quad
 \tilde{B} \,=\, \frac{e\, B}{m_e \,\omega_{pe}} \,=\,
 \frac{\omega_{ce}}{\omega_{pe}}, \quad \tilde{f_e} \,=\, f_e \,\frac{
 v_{th}^3}{ n_0}.
\end{equation}
\end{itemize}
This leads to the following dimensionless equations
\begin{equation}
\label{eq1}
\left\{
\begin{array}{l}
\displaystyle{\frac{\partial f_e}{\partial t} \,+\, \nabla_{{\bf x}}
\cdot \left ( {\bf v} f_e \right ) - \nabla_{\bf v} \cdot \left (({\bf
E}+ {\bf v} \times {\bf B} )f_e \right ) \,=\, \frac{\nu}{Z}
\,C_{e,e}(f_e,f_e) \,+\, \nu \,C_{e,i}(f_e)}, \\ \; \\
\displaystyle{\frac{\partial {\bf E}}{\partial
t}\,-\,\frac{1}{{\beta}^2}\, \nabla_{\bf x} \times {\bf B} = n\,{\bf
u},} \\ \; \\ \displaystyle{\frac{\partial {\bf B}}{\partial t} \,+\,
\nabla_{\bf x} \times {\bf E} \,=\, 0,} \\ \, \\
\displaystyle{\nabla_{\bf x} \cdot {\bf E} \,=\, (1\,-\,n), \quad
\nabla_{\bf x} \cdot {\bf B} \,=\, 0,}
\end{array}\right.
\end{equation}
where $\beta \,=\,v_{th} / c$, $\nu$ is the ratio between electron-ion collision frequency and plasma frequency
$$ \nu \,=\, \frac{Z\, n_0\, e^4 \ln\Lambda}{8\, \pi\,\epsilon_0^2
\,m_e^2\, v_{th}^3\,\omega_{pe}} \,=\, \frac{Z\,\ln\Lambda}{8\,
\pi\, n_0\, \lambda_D^3} \,=\, \frac{\nu_{e,i}}{\omega_{pe}} \quad{\rm
with }\quad \nu_{e,i} \,=\, \frac{Z\,n_0 \,e^4 \ln\Lambda}{8 \,\pi
\,\epsilon_0^2 \,m_e^2\, v_{th}^3}.
$$ The zero and first order moments of the distribution function are
$$ 
\left\{\begin{array}{l}
\displaystyle{n(t,{\bf x})\,=\, \int_{\RR^3} f_e(t,{\bf x},{\bf v}) {\,\bf dv},}
\\
\,
\\ 
\displaystyle{{\bf u}(t,{\bf x})\,=\, \frac{1}{n(t,{\bf x})} \int_{\RR^3} f_e(t,{\bf x},{\bf v}) \,{\bf
v}\, {\bf dv}.}
\end{array}\right.
$$ 
Moreover, in (\ref{eq1}) the dimensionless collision operators are considered
\begin{equation} 
  \label{eqc1}
\left\{\begin{array}{l}
 \displaystyle{C_{e,e}(f_e,f_e) \,=\, \nabla_{\bf v} \cdot \left(\int_{\RR^3}
 \Phi({\bf v}-{\bf v'}) \left [ f_e({\bf v}') \nabla_{{\bf v}}
 f_e({\bf v}) - f_e( {\bf v}) \nabla_{\bf v'} f_e({\bf v}') \right ]
 {\bf d v'} \right),}
\\
\,
\\
\displaystyle{ C_{e,i}(f_e) \,=\, \nabla_{\bf v} \cdot \left [ \Phi({\bf v})
  \nabla_{\bf v} f_e({\bf v}) \right],}
\end{array}\right.
\end{equation}
with $\Phi$ given by (\ref{op}).

\subsection{Test 1 : $1D$ two-stream instability}

The ICF physics involves a propagation of electron beams in plasma.
The plasma response to the beam consists in a return current that goes
opposite to the beam in order to preserve the quasineutrality. This
leads to a very unstable configuration favorable to the excitation of
plasma waves. We focus here on the instability with a perturbation
wavevector parallel to the beam propagation direction, namely the
two-stream intability.  Of course, this stands as an academic test
case but it is closely related to the physics of the ICF.  Also it is
a very demanding test for numerical schemes of transport, that have to
be specially designed (see Proposition~\ref{prop:3.2}). In particular,
a discrete dispersion relation relative to that problem is developed
to justify numerical choices for the second order scheme. For this
scheme also, during the limitation procedure, an additional
dissipation at extrema is introduced, compared to \cite{crouseilles}, in
order to preserve the solution from spurious oscillations. We will
show the sensitivity of the scheme with respect to the chosen limiter,
for this particular test case. Moreover, the fourth order scheme is
introduced to reduce numerical heating, for simulations intended to
deal with the two-stream instability.

The ($1D_{x}\times 1D_{v}$) Vlasov-Amp\`ere system
\eqref{simple2stream} is approximated on a Cartesian grid. For this
test case, we consider the scaling
(\ref{scaling})-(\ref{scaling1}). The initial distribution function
and electric field are
$$
\left\{
\begin{array}{l} 
\displaystyle{f^0({x},{v}) \,=\, \frac{1}{2} \left[
(1+A\,\cos(k{x}))\mathcal{M}_{1,v_d}({v})\,+\,
(1-A\cos(k{x}))\mathcal{M}_{1,-v_d}({v}) \right], } \\ \, \\
\displaystyle{E^0({x}) \,=\,0,}
\end{array}\right.
$$
where 
$$ \mathcal{M}_{1,v_d}(v)\,=\,\frac{1}{\sqrt{2\pi}}\,{\rm
e}^{-\|v-v_d\|^2/2}
$$
is the Maxwellian distribution function centered around $v_d$.

In order to compare the numerical heating associated with the second order and the fourth order scheme, we choose a strong perturbation amplitude $A\,=\,0.1$. The perturbation wavelength is
$k\,=\,2\pi / L$ and the beam initial mean velocities are $v_d\,=\,
\pm 4$, $L\,=\,25$ being the size of the periodic space domain. We
choose a truncation of the velocity space to be $v_{max}\,=\,12$ and
time steps are chosen to be $\Delta t = 1/200$.

The objectives of this numerical simulation are on the one hand to
compare the second order finite volume scheme (specially designed to
conserve exactly the discrete total energy, exept if the slope
limiters are active) for different slope limiters and the fourth order
MUSCL scheme. On the other hand we want to explore the effect of a
reduced number of grid points on the discrete invariants conservation.

\begin{figure}[htbp]
\begin{tabular}{cc}
\includegraphics[width=7.25cm,height=6.cm,angle=0]{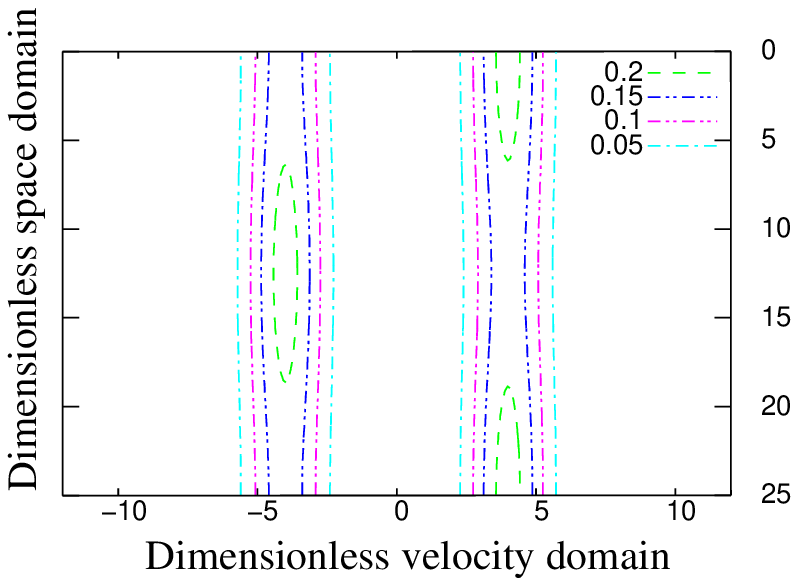}
&
\includegraphics[width=7.25cm,height=6.cm,angle=0]{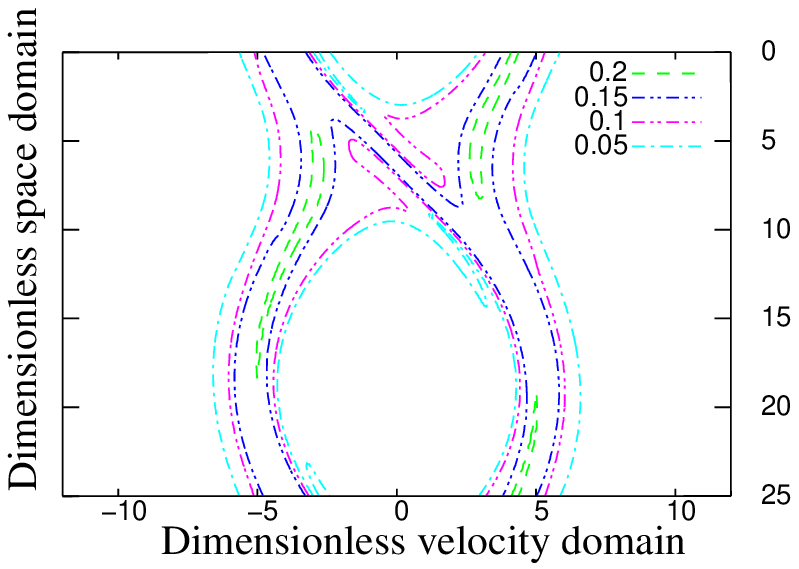}
\\
$(a)$ & $(b)$
\end{tabular}

\caption{{\it Beams phase space  $(a)$ at initial time, $(b)$ at $20$ plasma
periods (after saturation)}}
\label{structure}
\end{figure}

In Figure~\ref{structure}, two countersteaming beams that are
initially well separated in the phase space $(a)$ start to mix
together. They finally create a complicated vortex structure,
involving wave-particle interactions. This behaviour remains
quantitatively the same whatever the transport scheme is (second or
fourth order). However with a reduced number of grid points (smaller
than $128$ points in velocity), the second order (with limiter
\eqref{slopedef+}) and fourth order schemes present a different
behaviour for the total electric energy and total energy.

\begin{figure}[htb]
 \begin{tabular}{cc}      
\includegraphics[angle=0,height=6cm,width=7.25cm]{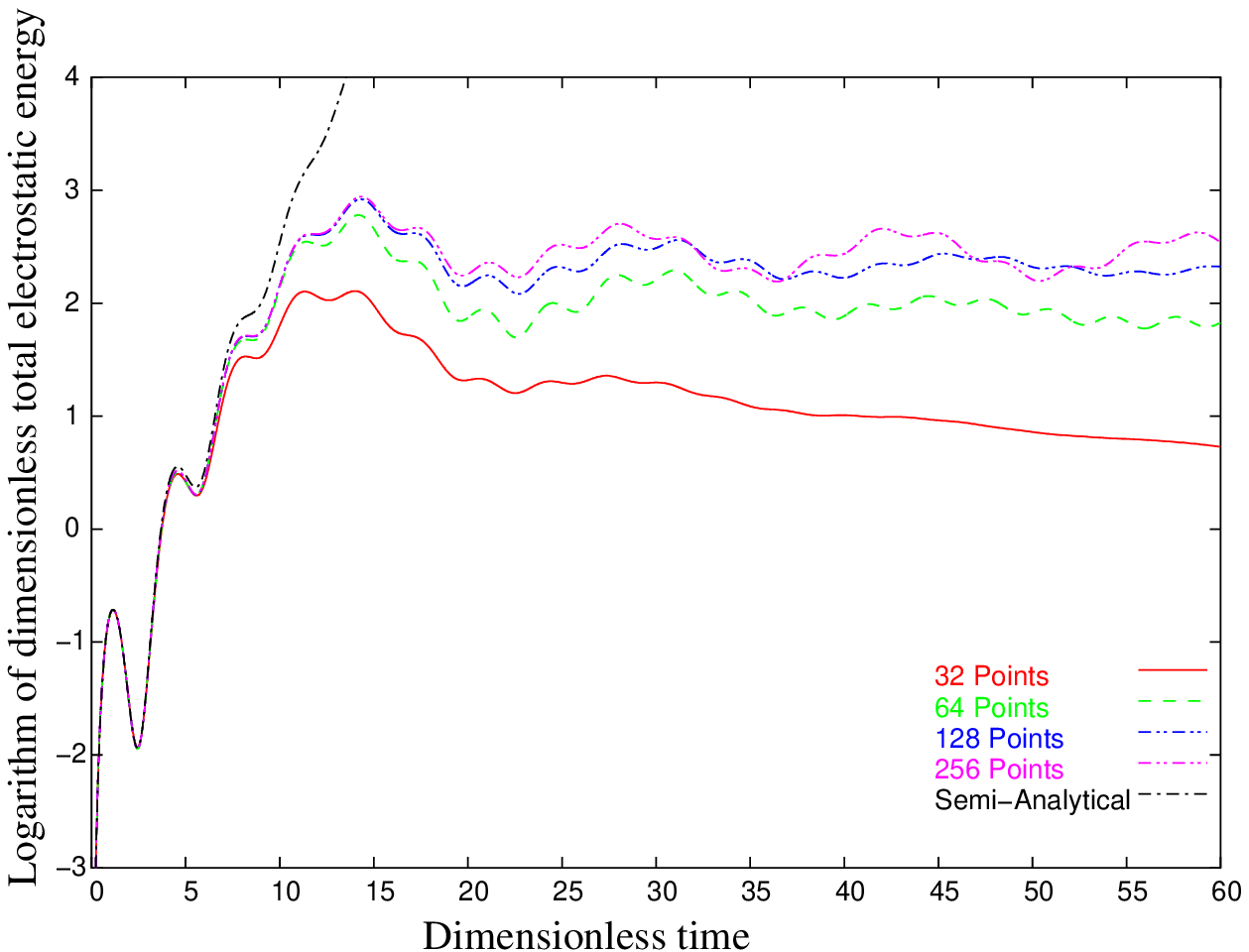}&
\includegraphics[angle=0,height=6cm,width=7.25cm]{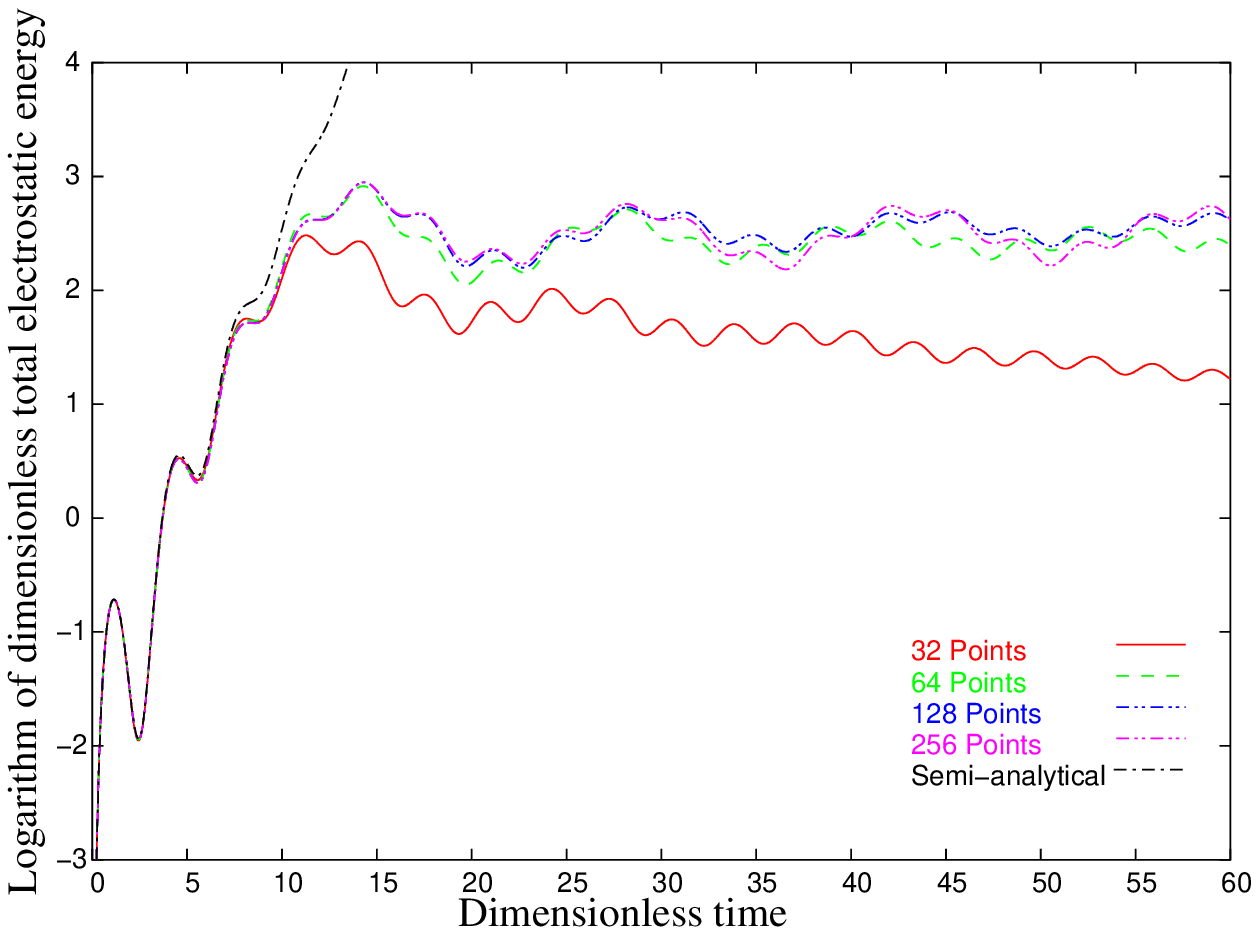}
\\
$(a)$ & $(b)$
\end{tabular}
\caption{{\it Evolution of the electrostatic energy for $32^2$,
$64^2$, $128^2$, $256^2$ grid points, and the semi-analytical solution
in the linear regime. Results are shown for $(a)$ the second order
with limiter \eqref{slopedef+} and $(b)$ fourth order transport
scheme }}
\label{field1}
\end{figure}

\begin{figure}[htbp]
\begin{tabular}{c}
\includegraphics[width=6.cm,height=7.25cm,angle=-90]{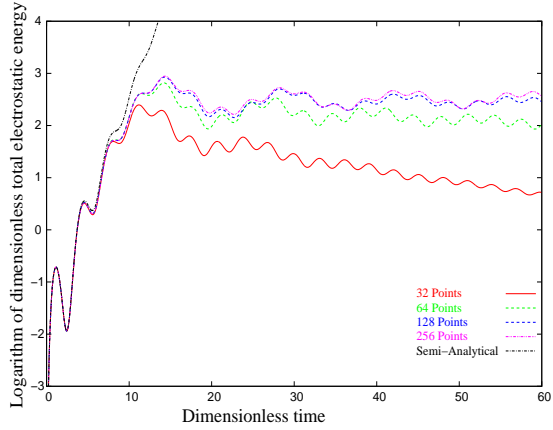} 
\end{tabular}
\caption{\it Evolution of the electrostatic energy for $32^2$, $64^2$,
$128^2$, $256^2$ grid points, and the semi-analytical solution in the
linear regime. Results are shown here for the second order scheme with
the limiters \eqref{slopedef+1}, with $b=2$.}
\label{fieldnewlimiter}
\end{figure}

\begin{figure}[htbp]
\begin{tabular}{cc}
\includegraphics[width=7.25cm,height=6.cm,angle=0]{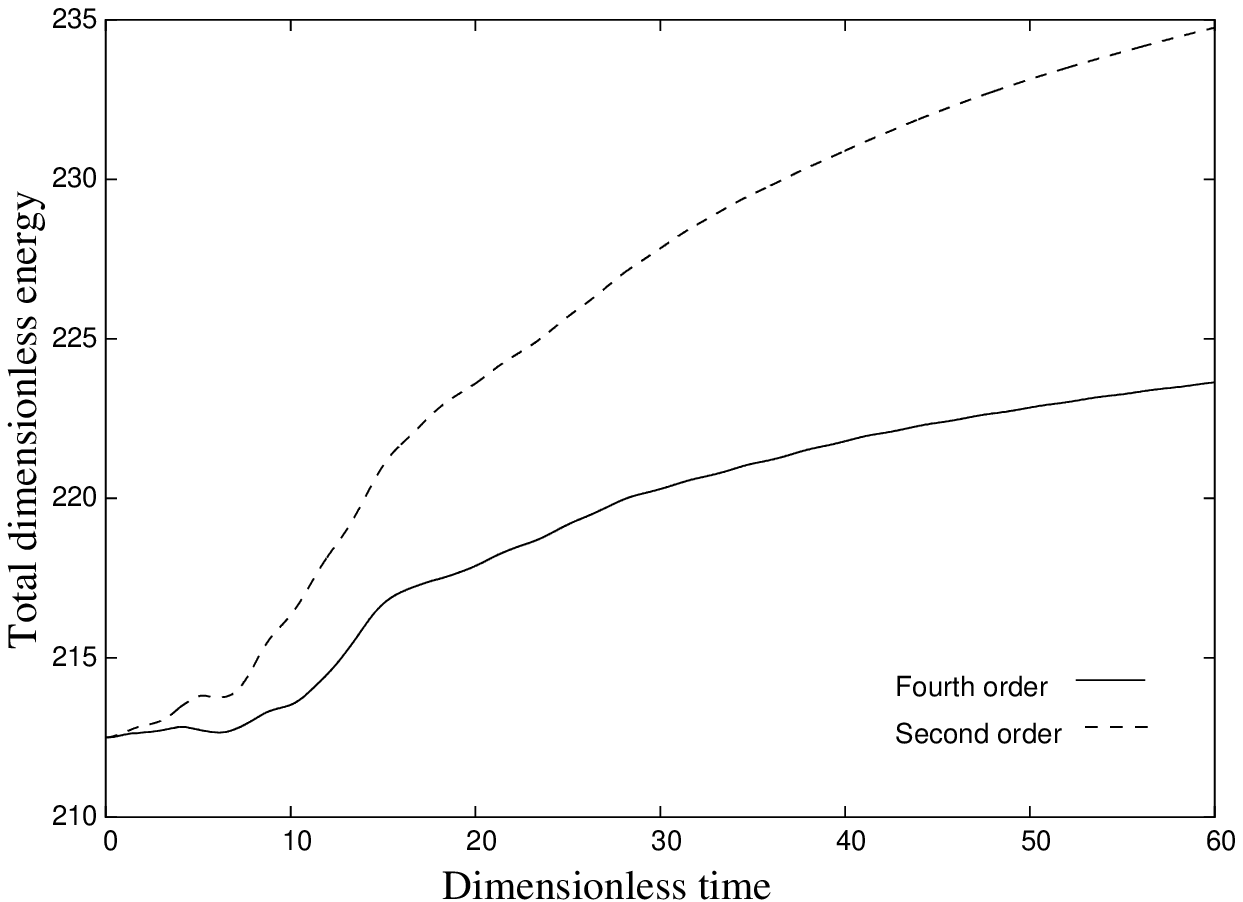}
&
\includegraphics[width=7.25cm,height=6.cm,angle=0]{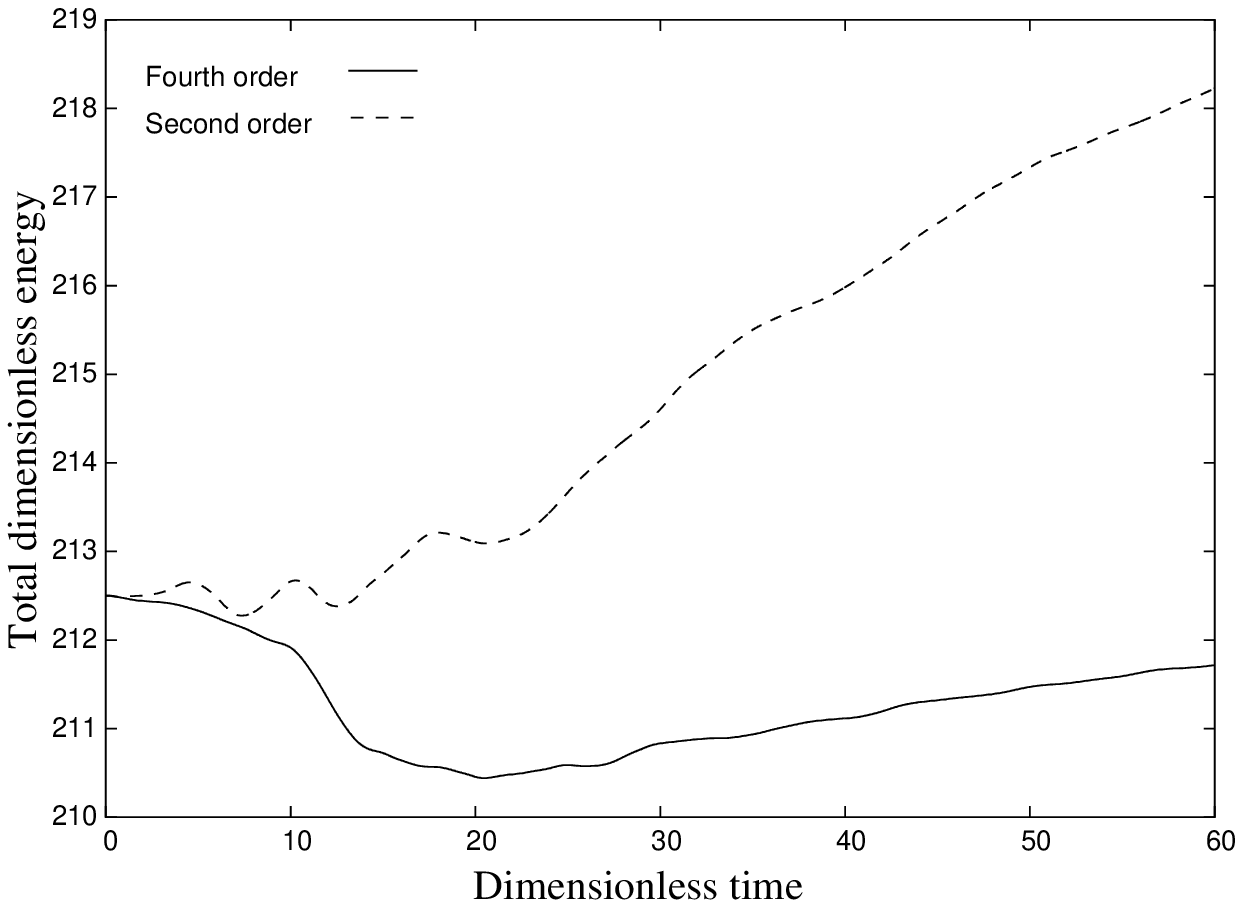}
\\
$(a)$  & $(b)$ 
\end{tabular}
\caption{\it Comparison of the energy evolution for the second (with
limiter \eqref{slopedef+}) and fourth order transport schemes. Results
are shown $(a)$ for $32^2$ $(b)$ $64^2$ grid points }
\label{energy}
\end{figure}

\begin{figure}[htbp]
\begin{tabular}{cc}
\includegraphics[width=6.cm,height=7.25cm,angle=-90]{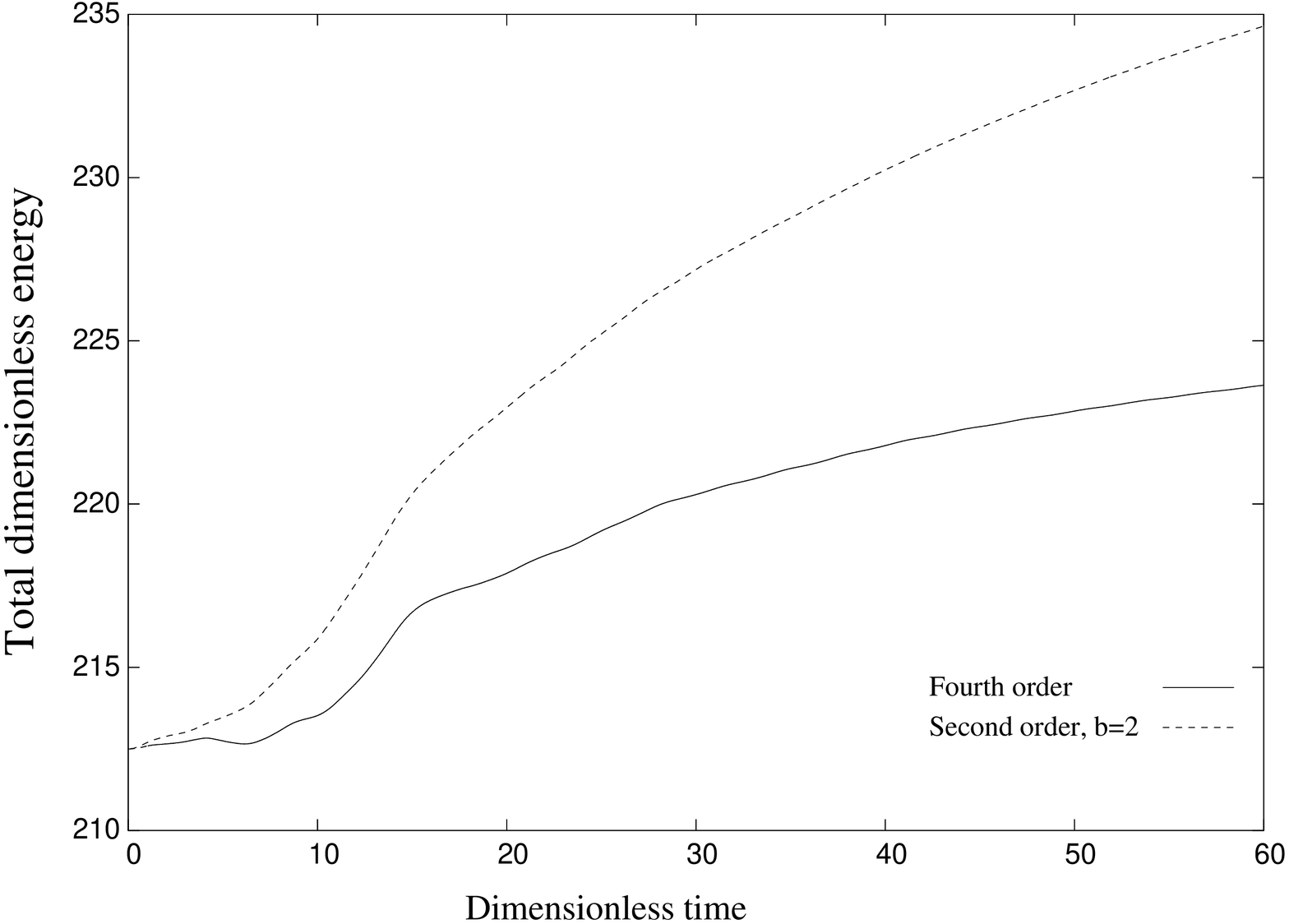}
&
\includegraphics[width=6.cm,height=7.25cm,angle=-90]{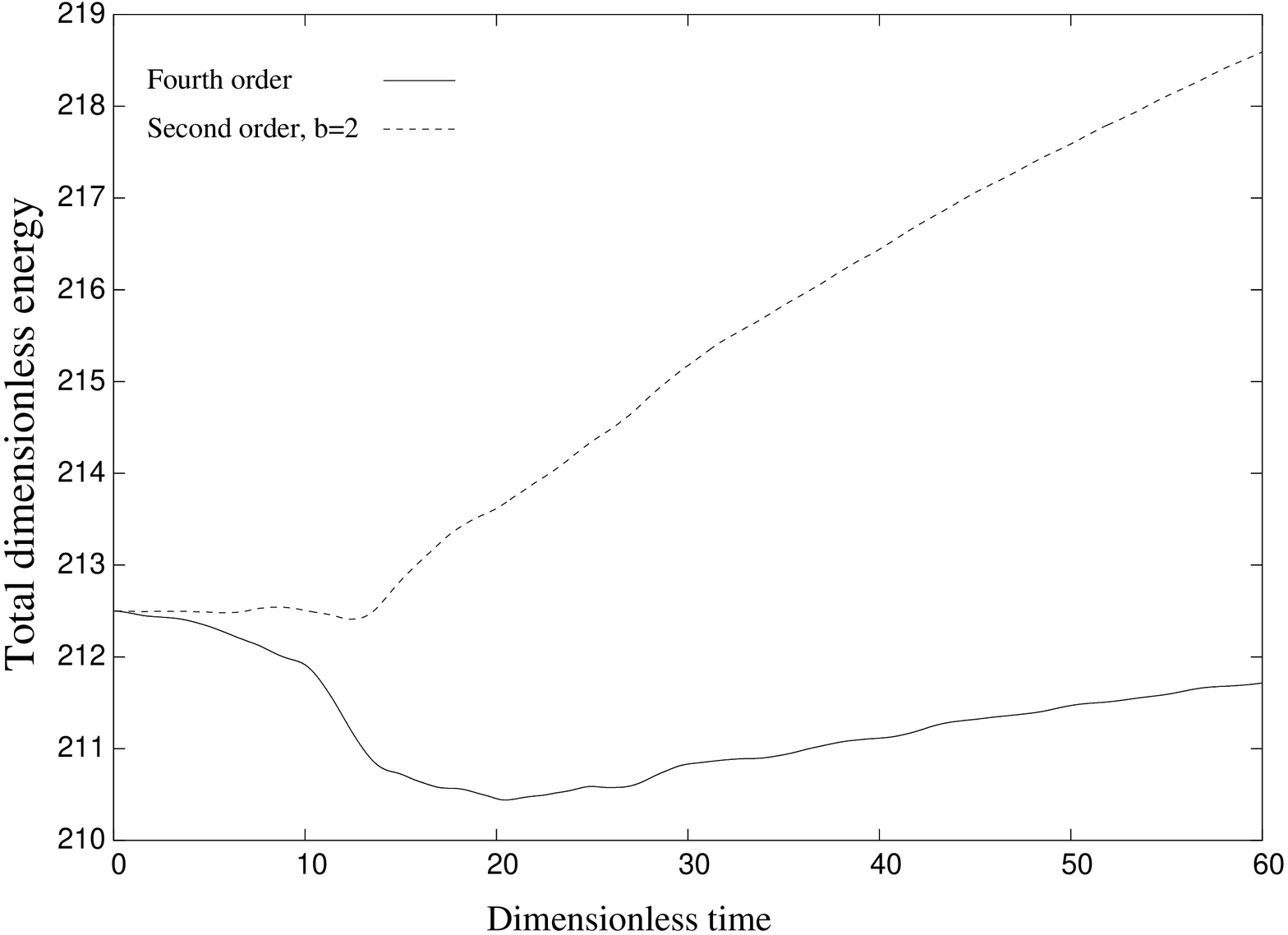}
\\
$(a)$  & $(b)$ 
\end{tabular}
\caption{\it Comparison of the energy evolution for the second (with
limiter \eqref{slopedef+1}, $b=2$) and fourth order transport
schemes. Results are shown $(a)$ for $32^2$ $(b)$ $64^2$ grid points }
\label{energy2}
\end{figure}

For reduced grid resolutions, of $32^2$ or $64^2$ points, the fourth
order scheme proves to be better than the second order one. For $32^2$
points, plasma oscillations at the plasma frequency in the nonlinear
phase are not reproduced with the second order scheme whereas they can
be seen with the fourth order scheme (see Figure~\ref{field1}).
Moreover for this resolution, the transition from the linear phase to
the nonlinear phase occurs earlier than it should for the second order
scheme.\\ As the grid resolution increases, the accuracy remains
better for the fourth order scheme than for the second order one in
the nonlinear phase (Figure \ref{field1}). The convergence toward
curves with $128^2$ or $256^2$ resolution grid is indeed better.  We
recall that quantities in Figure~\ref{field1} and \ref{fieldnewlimiter}
are plotted with a logarithmic scale, that smoothes out discrepancies
between curves.  In addition to these results, the respect of total
discrete energy conservation proves to be better for the fourth order
scheme than for the second order one at a reduced grid resolution, see
Figure~\ref{energy} and \ref{energy2}.

The use of limiters \eqref{slopedef+1} for the second order scheme
introduces accuracy improvements on the convergence behaviour and
capture of plasma wave structure at reduce grid resolutions, see
Figure~\ref{fieldnewlimiter}. However, the energy dissipation remains
quantitatively the same as the second order scheme with limiter
\eqref{slopedef+}, see Figures~\ref{energy} and \ref{energy2}.

  As this test case requires both a good preservation of invariants
and accuracy when nonlinear phenomena occur, we might conclude that
the fourth order scheme, with a resolution along each velocity
direction greater than $32$ cell, is well suited for our physical
applications.  The semi-analytical solution in the linear regime shown
in Figure~\ref{field1}, using a Green function, brings some
improvements compared to the classical validation in the linear
regime, based on instabilities growth rates in the linear regime. In
particular it discriminates precisely in time the linear and nonlinear
phases.

\subsection{Test 2: $1D$ X-mode plasma in a magnetic field}

This test case stands as a validation in the linear regime for the
coupling between Vlasov and Maxwell equations without collisions.  A
particular initial data is chosen (see the derivation in the
appendix~\ref{annexeC}) to trigger an X-mode plasma wave at a
well-defined frequency $\omega$.  This type of wave presents a mixed
polarization (longitudinal and transverse with respect to the magnetic
field), that propagates in the plane $P_{\perp}$, perpendicular to the
magnetic field direction.

The chosen frequency $\omega$ is a solution of the dispersion relation
\eqref{testtest1} of the linearized Vlasov-Maxwell equations,
introducing the equilibrium state $f^{(0)}\left ( \| {\bf v}
\|^2\right ) $.  The initial data are chosen such that $f^{(0)}$,
$\hat E_1$, $\hat E_2$, and $\hat B_3$ only depend on $\omega$,
$B^{(0)}$, $k_1\,=\,2\pi/L_1$ and $A$; where $ \hat f_n $, $\hat B_3
$, $\hat E_1 $ and $\hat E_2 $ are the reconstructed (in the
appendix~\ref{annexeC}) Fourier transforms of the distribution
function and electromagnetic fields. The magnetic field $B^{(0)}$ is
the nonperturbed magnitude of the magnetic field, $L_1$ is the length
of the space domain, $A$ is the perturbation amplitude. The initial
data can then be constructed with the help of truncated Fourier series
$$
\left \{
\begin{array}{l}
\displaystyle{ f^{(0)}({x_1},{\bf v}) \,=\, f^{(0)}(\|{\bf v}\|^2)
\,+\, \sum_{n=-2}^{2} \hat f_n({\bf v}_\perp) e^{i k_1 \,x_1 \,+\, i\,
n \psi}}, \quad x_1 \in (0,L_1), \, {\bf v}\in\RR^3, \\ \, \\
\displaystyle{ E_1(t,x_1) \, = \, \hat E_1 \,e^{-i \omega t + i k_1
x_1}}, \quad x_1 \in (0,L_1) \ , \\ \, \\ \displaystyle{ E_2(t,x_1) \,
= \, \hat E_2 e^{-i \omega t + i k_1 x_1}}, \quad x_1 \in (0,L_1) \ ,
\\ \, \\ \displaystyle{ B(t,x_1) \,=\; B^{(0)} \,+\, \hat B_{3} e^{- i
\omega t + i k_1 x_1}}, \quad x_1 \in (0,L_1).
\end{array}\right.
$$ We define $\psi$ as the angle in the cylindrical coordinates for
the velocity, defined with respect to the direction of the magnetic
field (See appendix \ref{annexeC}).
 
The normalisations are defined by relations
(\ref{scaling})-(\ref{scaling1}).  We choose $B^{(0)}=2$ and a rather
strong amplitude perturbation $A=0.1$ with periodic boundary
conditions on the space domain.  Also we have set
$\beta=v_{th}/c=0.05$.  The dispersion relation have been solved for
these parameters.  One of the solution $\omega \simeq 5.1432$ is
injected in the initial data set.

We considered $126$ points along the $1D$ space direction, and $64$
points along each velocity direction ${\bf v}=\,^t(v_1,v_2,v_3)$.  The
dimension of the space domain is $L_1\,=\,25$ whereas the truncation of
the velocity space occurs at $v_{max}=7$ for each velocity direction.
Furthermore, the time step is $\Delta t =1/200$.

\begin{figure}[htb]
\centering
\includegraphics[height=12cm,width=12cm,angle=-90]{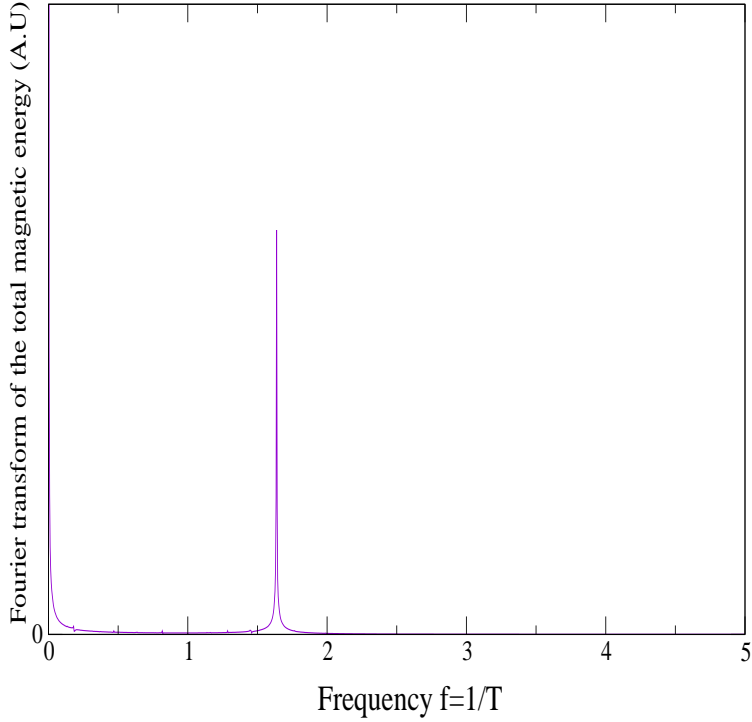}
\caption[]{{\it Discrete Fourier spectrum in frequency of the discrete
analogous of the total dimensionless magnetic energy $ \displaystyle
\int_{0}^{L_1} \frac{\|B_3\|^2}{2} dx_1$.}}
\label{figfourier}
\end{figure}

\begin{figure}[htb]
 \centering
 \includegraphics[angle=0,width=14cm,keepaspectratio=true]{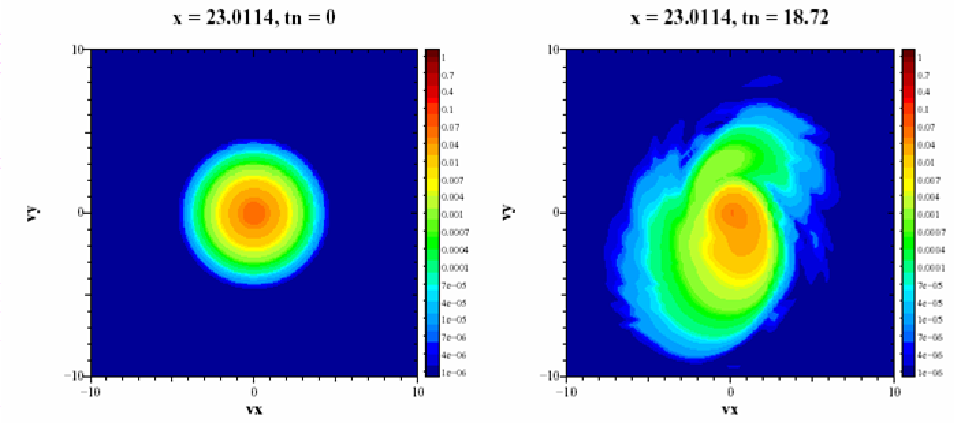}
\caption{\it Projection on the $v_1-v_2$ velocity domain for the
distribution function is shown at initial time $t_n=0$ and at
$t_n=18.72$, for a particular point of the space domain,
$x_1=23.0114$, $v_3=0$.}
\label{distribmag}
\end{figure}

The Fourier spectrum in Figure~\ref{figfourier} exhibits a well
defined frequency $f=1/T \simeq 1.6375$ (corresponding to a period
$T$) for the total magnetic energy, that corresponds to a frequency
$f/2$ for the magnetic field oscillations.  We finally find $\omega =
2 \pi\left ( \frac{f}{2}\right ) \simeq 5.1443$ from the numerical
solution, to be compared with the analytical results $5.1432$.  This
proves a good accuracy of the numerical results, while the
distribution function is greatly affected by the magnetic field.  As
an illustration, we show in Figure~\ref{distribmag} how the magnetic
field makes the distribution function rotate in the velocity space
perpendicular to the magnetic field axis.

\section{Approximation of the collision operators}
In the following, the presentation is restricted to the space homogeneous equation, for the sake of simplicity,
$$
\left\{
\begin{array}{l}
\displaystyle{\frac{\partial f}{\partial t} \,=\, C_{e,e}(f,f)
\,+\,C_{e,i}(f)}, \\ \, \\ f(0,{\bf v}) \,=\,f^{(0)}({\bf v}),
\end{array}\right.
$$  
where $C_{e,e}(f,f)$ and  $C_{e,i}(f)$ are  given by (\ref{eqc1}).

\subsection{Discretization of the Lorentz operator}
We consider $f_{\bf j}$ an approximation of the distribution function
$f{\bf( v_j)} $ and introduce the operator ${\bf D}$, which denotes a
discrete form of the usual gradient operator $\nabla_{\bf v}$ whereas
${\bf D^*}$ represents its formal adjoint, which represents an
approximation of $-\nabla_{\bf v} \cdot$. Therefore, for any test
sequence $(\psi_{\bf j})_{{\bf j} \in \mathbb{Z}^3}$, we set $({\bf
D}\psi_{\bf j})_{{\bf j} \in \mathbb{Z}^3} $ as a sequence of vectors
of $\RR^3$
$$ {\bf D}\psi_{\bf j} \,=\, {}^{t}(D_1\psi_{\bf j}, D_2\psi_{\bf j},
D_3\psi_{\bf j}) \in \RR^3,
$$ where $D_s$ is an approximation of the partial derivative
$\frac{\partial}{\partial v_{s}}$ with $s \in \{1,2,3\}$. In order to
preserve the property of decreasing entropy at the discrete level, we
use the log weak formulation of the Lorentz operator \cite{degond}
$$ 
\int_{\RR^3}\,C_{e,i}(f)({\bf v})\,\psi({\bf v})\,{\bf dv} \,=\,
-\int_{\RR^3}\Phi({\bf v})\,f({\bf v})\,\nabla_{\bf v} \log(f({\bf
v})) \cdot \nabla_{\bf v} \psi({\bf v})\,{\bf dv},
$$ where $\Phi$ is given by (\ref{op}) and $\psi$ is a smooth test
function. Then, using the notations previously introduced, the
discrete operator $C_{e,i}^{\Delta v}(f)$ is given by
\begin{eqnarray}
\label{lorentzlogweak}
C_{e,i}^{\Delta v}(f)({\bf v_j}) \,=\, - {\bf D^*} \left[\frac{1}{\|{\bf
v}_{\bf j}\|^3}\, S({\bf \tilde{v}_j})\, f_{\bf j} \, {\bf
D}(\log(f_{\bf j}))\,\right],
\end{eqnarray} 
where $S(\tilde{{\bf v_j}}) $ is the following matrix
$$ 
S(\tilde{{\bf v_j}} ) \,=\, \|\tilde{{\bf v_j}}\|^2 \,{\rm Id}
\,-\, \tilde{{\bf v_j}} \otimes \tilde{{\bf v_j}}.
$$ Now, ${\bf \tilde{v}_j} $ has to satisfy the discrete conservation
of energy
\begin{eqnarray}
\label{compatdiscrete}
\frac{D_1(\|{\bf {v}_j}\|^2)}{{\bf \tilde{v}}_{j_1}} \,=\,
\frac{D_2(\|{\bf {v}_j}\|^2)}{{\bf \tilde{v}}_{j_2}} \,=\,
\frac{D_3(\|{\bf {v}_j}\|^2)}{{\bf \tilde{v}}_{j_3}}.
\end{eqnarray}
Then, we consider the $8$ uncentered operators ${\bf D^{{\bf
\epsilon}}}$, with the formalism:
$$ {\bf D^{ {\bf \epsilon} }}\,=\,{}^{t}(D^{{ \epsilon_1}}_1,D^{{
\epsilon_2}}_2,D^{{ \epsilon_3}}_3),
$$ with ${\bf \epsilon} ={}^{t}(\epsilon_{1},\epsilon_{2},\epsilon_{3})$,
and $\epsilon_{i} \in \{+1,-1\}$ for $i \in \{1,2,3\} $. More
precisely, the operator $ D^{\epsilon_{i}} $ is the forward uncentered
discrete operator if $\epsilon_{i}=+1$ and the backward uncentered
discrete operator if $\epsilon_{i}=-1$:
\begin{eqnarray}
\label{forwardbackward}
{\bf D^{\epsilon}} {\bf \Psi}_{{\bf j}} = \frac{1}{\Delta v} \left (
\begin{array}{ccc}
\epsilon_{1}[\Psi_{j_1+\epsilon_{1}}-\Psi_{j_1}] \\
\epsilon_{2}[\Psi_{j_2+\epsilon_{2}}-\Psi_{j_2}] \\
\epsilon_{3}[\Psi_{j_3+\epsilon_{3}}-\Psi_{j_3}]
\end{array} \right )
\end{eqnarray}
This $8$ operators respectively match to $8$ expressions of ${\bf
\tilde{v}_j^{\epsilon}}$, following (\ref{compatdiscrete})
$$  
{\bf \tilde{v}_j^{\epsilon}} \,=\, \frac{1}{2} \left (
 {\bf{v}_j} \,+\,  {\bf{v}_{j+{\bf \epsilon}}} \right ).          
$$ This choice has been made to avoid the use of the centered discrete
operator that conserves non physical quantities.  On the other hand,
the uncentered operators, taken separately, introduce some artificial
unsymmetry in the distribution function leading to a loss of accuracy
when coupling to Maxwell equations.  To overcome these difficulties,
following the idea of \cite{buetcordier}, we introduce a
symmetrization of the discrete operator based on the averaging over
the eight uncentered discretizations:
\begin{eqnarray}
\label{lorentzlogweakavrg}
& \displaystyle C^{\Delta v}_{e,i}(f)({\bf v_j}) \,=\, \frac{1}{8}
    \sum_{{\bf \epsilon}} C_{e,i}^{\epsilon}(f) \nonumber \\ &
    \displaystyle C_{e,i}^{\epsilon}(f) = -{\bf D^{*,{\bf \epsilon}}}
    \left[ \frac{1}{\|{\bf v_j}\|^3}S({\bf \tilde{v}_j^{{\bf
    \epsilon}}}) \,f_{\bf j} \,{\bf D^{{\bf \epsilon}}}(\log(f_{\bf
    j})) \right]. \nonumber
\end{eqnarray} 
This final expression will introduce an additional discrete symmetry
property compared to the operator presented in \cite{crouseilles}.

 We now present the discrete properties for the electron-ion collision
operator. We have the classical properties: mass and energy
preservation, an entropy decreasing behaviour, the positivity
preservation of the distribution function in a finite time sequence.
The proofs are not detailed here but can be deduced easily from those
presented in \cite{crouseilles}.  The difference stands in the fact
that we obtain the operator as an average over the full set of the
uncentered operators (instead of an average over two operators). This
modification allows to get a discrete analogous of the symmetry
property presented in Proposition~\ref{prop:01}:

\begin{proposition}
 Under the condition (\ref{compatdiscrete}) on $ {\bf \tilde{v}_j}$,
 the discretization (\ref{lorentzlogweakavrg}) to the Lorentz operator
 (\ref{Cei}) satisfies the following properties,
\begin{itemize}
\item it preserves mass and energy,
\item it decreases discrete entropy
$$ H(t)\,=\,\Delta v^3 \,\sum_{{\bf j} \in \mathbb{Z}^3} f_{{\bf
j}}(t)\,\log(f_{{\bf j}}(t)),
$$
\item there exists a time-sequence $\Delta t_n $ such that the scheme
$$ f_{{\bf j}}^{n+1} \,=\, f_{{\bf j}}^{n} \,+\, \Delta t \,C^{\Delta
v}_{e,i}(f)({\bf v_j}),
$$ defines a positive solution at any time {\it i.e.} $\sum_n t_n
\,=\,+\infty$.
\end{itemize}

Furthermore, if $f_{\bf j}$ is symmetric with respect to $0$ in the
direction $j_k$ at time $t^n$, then this property is preserved at time
$t^{n+1}$,
\begin{equation}
\label{symmetrypropdiscrete}
\sum_{{\bf j}\in \ZZ^3} C^{\Delta v}_{e,i}(f)({\bf v_j})v_{j_k} \Delta v^3
\,=\,0.
\end{equation} 
\end{proposition}

\begproof We prove the last property and rewrite the operator
(\ref{lorentzlogweakavrg}) in a different manner, assuming we have a
symmetry along the velocity direction $v_{j_k}$
\begin{equation}
\label{lorentzredistrib}
C^{\Delta v}_{e,i}(f)({\bf v_j})\,=\,\frac{1}{8}\sum_{\epsilon }
C_{e,i}^{{\bf \epsilon}}(f)({\bf v_j})\,=\,\frac{1}{4} \left [
\sum_{} \frac{1}{2}
\left ( C_{e,i}^{{\bf \epsilon}^{+,(k)}}({\bf v_j}) + C_{e,i}^{{\bf
{\epsilon}}^{-,(k)}}({\bf v_j}) \right ) \right ] ,
\end{equation}
where the notation ${\bf \epsilon}^{\pm,(k)}$ refers to
\begin{eqnarray}
\label{epsilondef}
\,\left \{ \begin{array}{ll} {\bf \epsilon}^{\pm,(k)}_{i} = \pm 1 &
\textrm{ if } i = k, \\ {\bf \epsilon}^{\pm,(k)}_{i} = {\bf
\epsilon}_{i} & \textrm{ if } i \ne k.
 \end{array} \right .
\end{eqnarray}
We are interested in the cancellation of the operator $\displaystyle
\sum_{{\bf j}\in \ZZ^3} C^{\Delta v}_{e,i}(f)({\bf v_j})v_{j_k}$. This
is equivalent to the cancellation of
\begin{eqnarray*}
  Q^{{(k)}} &:=& \sum_{{\bf j} \in \mathbb{Z}^3} \left ( C_{e,i}^{{\bf
    \epsilon}^{+,(k)}}({\bf v_j}) + C_{e,i}^{{\bf {\epsilon}}^{-,(k)}}
    ({\bf v_j}) \right ) v_{j_k}  \\ &=& \sum_{{\bf j} \in \ZZ^3}
    \frac{1}{\|{\bf v_j}\|^3} f_{\bf j} \left [ S(\tilde{\bf
    v_j}^{{\bf \epsilon}^{+,(k)}}) {\bf D}^{{\bf
    \epsilon}^{+,(k)}}\log(f_{\bf j})\right ]\cdot {\bf D}^{{\bf
    \epsilon}^{+,(k)}} v_{j_k} \\ &+& \sum_{{\bf j} \in \ZZ^3}
    \frac{1}{\|{\bf v_j}\|^3} f_{\bf j} \left [ S(\tilde{\bf
    v_j}^{{\bf \epsilon}^{-,(k)}}) {\bf D}^{{\bf
    \epsilon}^{-,(k)}}\log(f_{\bf j}) \right ]\cdot {\bf D}^{{\bf
    \epsilon}^{-,(k)}} v_{j_k}.
\end{eqnarray*}
Then, since $ {\bf D}^{{\bf \epsilon}^{+,(k)}} v_{j_k} = {\bf
D}^{{\bf \epsilon}^{-,(k)}} v_{j_k} ={\bf e_k}$, it yields
\begin{eqnarray*}
Q^{{(k)}} &=& \sum_{{\bf j}\in\ZZ^3} \frac{1}{\|{\bf v_j}\|^3} f_{\bf
  j} \left ( \sum_{i\ne k} \left (
  \tilde{v}_{j_i}^{{\epsilon}^{+,(k)}_{i}} \right )^2 \right
  )D^{\epsilon_{k}^{+,(k)}}(\log(f_{\bf j})) \\ &-& \sum_{{\bf j} \in
  \mathbb{Z}^3}\frac{1}{\|{\bf v_j}\|^3} f_{\bf j}
  \tilde{v}_{j_k}^{\epsilon_{k}^{+,(k)}} \left ( \sum_{i\ne k}
  \tilde{v}_{j_i}^{\epsilon_{i}^{+,(k)}} D^{\epsilon_{i}^{+,(k)}}(\log
  f_{\bf j}) \right ) \\ &-& \sum_{{\bf j} \in
  \mathbb{Z}^3}\frac{1}{\|{\bf v_j}\|^3} f_{\bf j} \left ( \sum_{i\ne
  k} \left ( \tilde{v}_{j_i}^{{\epsilon}^{-,(k)}_{i}} \right )^2
  \right )D^{{\epsilon}^{-,(k)}_{k}}(\log(f_{\bf j})) \\ &-&
  \sum_{{\bf j} \in \mathbb{Z}^3}\frac{1}{\|{\bf v_j}\|^3} f_{\bf j}
  \tilde{v}_{j_k}^{{\epsilon}^{-,(k)}_{k}} \left ( \sum_{i\ne k}
  \tilde{v}_{j_i}^{{\epsilon}^{-,(k)}_{i}}
  D^{{\epsilon}^{-,(k)}_{i}}(\log f_{\bf j}) \right ).
\end{eqnarray*}
Then using definition (\ref{epsilondef}) and the symmetry of $f_{\bf
j}^n$ with respect to $0$ in the velocity direction $v_{j_k}$, we
obtain $ Q^{{(k)}}=0$. Then multiplying (\ref{lorentzredistrib}) by $
v_{j_k}$ and integrating in the full velocity space gives the relation
(\ref{symmetrypropdiscrete}). This relation implies that $f_{\bf
j}^{n+1}$ is symmetric with respect to 0 in the direction $v_{j_k}$.
\endproof

\subsection{Discrete Landau operator}

We consider the discretization of the FPL operator (\ref{Cee}) on the
whole 3D velocity space. It is based on the entropy conservative
discretization introduced in \cite{degond}, where a discrete weak log
form of the FPL operator is used. This discretization yields:

\begin{equation}
\label{weaklog}
\left\{
\begin{array}{l}
\displaystyle{\frac{df_{\bf j}(t)}{dt}= \left ( {\bf D^*}\rho(t)\right
)_{\bf j} \qquad {\bf j} \in \ZZ^3,} \\ \, \\ \displaystyle{ \rho(t) =
\Delta v^3 \sum_{{\bf m} \in \ZZ^3} f_{\bf j}(t)f_{\bf m}(t)\Phi({\bf
v_j}-{\bf v_m}) \left ( {\bf D}(\log(f(t))_{\bf j} - {\bf D}(\log
f(t))_{\bf m} \right),}
\end{array} \right.
\end{equation}
where ${\bf D}$ stands for a downwind or upwind finite discrete
operator approximating the usual gradient operator $\nabla_{\bf
v}$. This uncentered approximation ensures that the only equilibrium
states are the discrete Maxwellian. The use of centered discrete
operators would have lead to non physical conserved quantities. The
discretization of the FPL operator is then obtained as the average
over uncentered operators, but here for a different reason as in the
previous section, on the electron-ion collision operator
discretization.  In \cite{buetcordier1}, the scheme is rewritten as
the sum of two terms: a second order approximation and an artificial
viscosity term in $\Delta v^2$ which kills spurious
oscillations. However the computational cost of a direct approximation
of (\ref{weaklog}) remained too high. Therefore, a multigrid technique
has been used. We refer to \cite{buetcordier1} and \cite{buetcordier2}
for the details of the implementation on the FPL operator.
Nevertheless, these latter approaches introduce a new approximation
than can affect accuracy.  Based on \cite{Pareshi}, Crouseilles and
Filbet proposed another approach and noticed that the discrete FPL
operator (\ref{weaklog}) in the Fourier space can be written as a
discrete convolution, which directly gives a fast algorithm.
Here we adopt the multigrid method, detailed in \cite{buetcordier1},
that has a complexity of order ${O}(n_v^3\log n_v^3)$.

This discrete approximation preserves positivity, mass, momentum, energy, and decreases the entropy. Moreover the discrete equilibrium states are the discrete Maxwellian.

\section{Numerical results}

\subsection{Scaling with collision frequency}
 For the analysis of collisional processes, a new scaling is
introduced here, that allows time steps to be of the order of the
electron-ion collision time. In order to account for transport
phenomena occuring at the collision time scale, several parameters are
required: the electron-ion collision frequency $\nu_{e,i} $, the
associated mean free path $\lambda_{e,i}$, the thermal velocity $
v_{th}$, and the cylotron frequency $\omega_{ce} $
\begin{equation} 
  \label{eq_def_n2}
  \nu_{e,i} = \frac{Z\, n_0 \,e^4 \,\ln\Lambda}{8 \,\pi \,\epsilon_0^2\, m_e^2 \,v_{th}^3}
  \, , \quad 
  \lambda_{e,i} = \frac{v_{th}}{ \nu_{e,i}} \, , 
  \quad v_{th} = \sqrt{\frac{\kappa_B T_0}{m_e}}\,  , 
  \quad \omega_{ce} = {\frac{e B}{m_e}}\,  ,
\end{equation}
These parameters enable us to define the dimensionless parameters with tilde.
\begin{itemize}
\item Dimensionless time, space and velocity, respectively
\begin{equation}
 \label{scalingb0}
\tilde{t}=\nu_{e,i}t, \quad \tilde{x}= {x}/{\lambda_{e,i}}, \quad \tilde{v}={v}/{v_{th}}.
 \end{equation}
\item Dimensionless electric field, magnetic field, and distribution
function, respectively
\begin{equation}
  \label{scalingb1}
\tilde{E} = \frac{eE}{ m_e v_{th} \nu_{e,i}} \,, \quad \tilde{B} =
 \frac{e B}{m_e \nu_{e,i}} = \frac{\omega_{ce}}{\nu_{e,i}} \,, \quad
 \tilde{f_e} = f_e \frac{ v_{th}^3}{ n_0}.
\end{equation}
\end{itemize}
This leads to the following dimensionless equations
\begin{equation}
  \label{eqb1}
\left\{
\begin{array}{l}
\displaystyle{\frac{\partial f_e}{\partial t} + \nabla_{{\bf x}} \cdot \left (
{\bf v}f_e \right ) - \nabla_{\bf v} \cdot \left (({\bf E}+ {\bf v} \times
{\bf B} )f_e \right ) = \frac{1}{Z} C_{e,e}(f_e,f_e)+ C_{e,i}(f_e),}
\\
\,
\\
\displaystyle{\frac{\partial {\bf E}}{\partial t}-\frac{1}{{\beta}^2} \nabla_{\bf
x} \times {\bf B} = \frac{1}{{\alpha}^2}n{\bf u}, }
\\
\,
\\
\displaystyle{\frac{\partial {\bf B}}{\partial t}+ \nabla_{\bf x} \times {\bf E}
= 0,} 
\\
\,
\\
\displaystyle{ \nabla_{\bf x} \cdot {\bf E} = \frac{1}{{\alpha}^2}(1-n),}
\\
\,
\\
\displaystyle{\nabla_{\bf x} \cdot {\bf B} = 0,}
\end{array}\right.
\end{equation}
where $ \alpha=\nu_{e,i} / \omega_{pe}$ and $ \beta = v_{th} / c$. The collision terms $C_{e,e}(f_e,f_e)$ and $ C_{e,i}(f_e)$ are given in (\ref{eqc1}).

\subsection{1D temperature gradient test case}

\begin{figure}[htb]
       \centering
 \includegraphics[height=12cm,width=12cm,keepaspectratio=true]{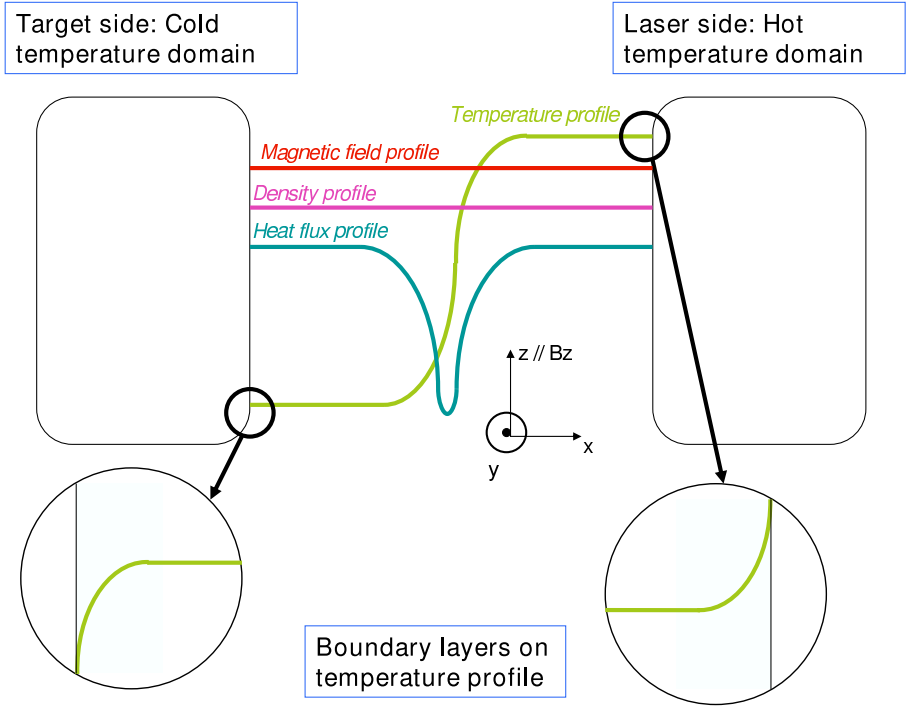}   
\caption[]{{\it Initial configuration for the temperature gradient
        test case: a temperature profile is considered between to two
        domains of plasma with particles at thermodynamical
        equilibrium. Zero current boundary conditions enable to
        maintain mass conservation. A heat flux is generated wherever
        there is a nonzero temperature gradient, as well as boundary
        layers on the heat flux, temperature, and electromagnetic
        profiles. \label{fig1} }}
\end{figure}

In the context of laser produced plasma, the heat conduction is the
 leading mecanism of energy transport between the laser energy
 absoption zone and the target ablation zone.\\ In such a system, the
 parameters of importance for the heat flux are
\begin{itemize}
\item The effective electron collision mean free path $ \lambda_e$.
\item The electron temperature gradient length $ \lambda_T$.
\item The magnetic field $B$ and its orientation with respect to
$\nabla T$.
\end{itemize}
These parameters enable to distinguish different regimes of transport,
according to the Knudsen and the Hall parameters.\\ The Knudsen number
$K_n$ is a mesure of the thermodynamical non-equilibrium of the
system

\begin{eqnarray} 
  \label{eqKN0}
&& K_n = \frac{\lambda_e}{\lambda_T}.
\end{eqnarray}
A regime characterized by $K_n \rightarrow 0$ refers to an
hydrodynamical descripion, whereas a regime characterized by $K_n \ge
1$ refers to a kinetic description, where the nonlocal phenomena
appear. The parameters for ICF imply $ K_n \ge 0.1 $, while the
classical, local approach fails at $ K_n \ge 0.01 $. This premature
failure of the classical diffusion approach in plasma is explained by
a specific dependence of the electron mean free path on their
energy. In our applications the energy is transported by the fastest
electrons, which have a much longer mean free path.\\ The Hall
parameter $ \chi= \omega_c \tau$ quantifies the relative importance of
magnetic and collisional effects. $\omega_c=eB/m_e$ is the electron
cyclotron frequency and $\tau$ the mean electron-ion collision time

\begin{eqnarray} 
 \label{dimbraginskii1}
&& \tau = \frac{3}{4}
\frac{16\pi^2\epsilon_0^2\sqrt{m_e}T_e^{3/2}}{\sqrt{2\pi}n_iZ^2e^4ln\Lambda}.
\end{eqnarray}
For this test case, a simple gradient temperature configuration is
shown in figure (\ref{fig1}), modelling the following situation:
through a layer of homogeneous plasma, a laser deposits its energy on
the hot temperature side and the absorbed energy is transported with
electrons to the cold temperature side.\\

Let us define the average over velocity of a function $A({\bf v})$

\begin{eqnarray} 
 \label{average}
<A> = \frac{1}{n_e} \int_{\RR^3} Af {\bf dv} \quad ,
\end{eqnarray}

where $n_e(t,{\bf x})=\displaystyle \int_{\RR^3} f(t,{\bf x},{\bf
v}){\bf dv}$ is the density of electrons.

Following \cite{braginskii,epperlein}, we introduce the
macroscopic quantities

\begin{eqnarray} 
 \label{macro}
\left\{
\begin{array}{l}
\displaystyle {\bf j} = q_en_e \left < {\bf v} \right >, \\ \, \\
 \displaystyle {\bf q} = \frac{1}{2}m_en_e \left < ({\bf v \cdot v })
 {\bf v}\right >,\\ \, \\ \displaystyle \label{transfertmom} {\bf R}=
 \int_{\RR^3} m_e{\bf v} C_{e,i}(f_e){\bf dv},
\end{array}\right.
\end{eqnarray}

\begin{eqnarray} 
 \label{macro1}
\left\{
\begin{array}{l}
\displaystyle p = n_e T_e = \frac{1}{3}m_en_e \left < ({\bf v-<{\bf
  v}>})\cdot({\bf v-<{\bf v}>}) \right >,\\ \, \\ \displaystyle {\bf
  \Pi} = \frac{1}{3}m_en_e \left < ({\bf v-<{\bf v}>})\otimes({\bf
  v-<{\bf v}>}) \right > -p{\bf I},\\ \, \\ \displaystyle {\bf
  q_{loc}} = \frac{1}{2}m_en_e \left < \left [ ({\bf v-<{\bf
  v}>})\cdot({\bf v-<{\bf v}>}) \right ] ({\bf v-<v>}) \right >.
\end{array}\right.
\end{eqnarray}
There, $ {\bf j}$ is the electric current, ${\bf q}$ the total heat
flow, ${\bf R}$ the friction force accounting for the transfer of
momentum from ions to electrons in collisions, $T_e$ is the
temperature, $p$ is the scalar intrinsic pressure, ${\bf \Pi}$ is the
stress tensor, ${\bf q_{loc}}$ is the intrinsic heat flow and ${\bf
I}$ the unit diagonal tensor.\\ Quantities $p$, ${\bf \Pi}$ and ${\bf
q_{loc}}$ are defined in the local reference frame of the electrons,
whereas ${\bf j}$, ${\bf q}$ and ${\bf R}$ are defined relative to the
ion center of mass frame. Ions are supposed to be at rest. We have the
relation

\begin{eqnarray}
 &&\label{trans} {\bf q_{loc}} = {\bf q} + {\bf j}\cdot
(\frac{5}{2}p{\bf I} + {\bf \Pi})/(n_ee) + {\bf
j}(\frac{1}{2}m_en_e<{\bf v}>^2)/(n_ee).
\end{eqnarray}

The validation of our Fokker-Planck solver in the domain close to the
hydrodynamical regime (local regime) requires knowledge of transport
coefficients. Following the formalism of Braginskii \cite{braginskii}
for the transport relations, the transport coefficients in the
hydrodynamical regime have been calculated by Epperlein in
\cite{epperlein}. These coefficients ${\bf \alpha_{ep}}$, ${\bf
\beta_{ep}}$, ${\bf \kappa_{ep}}$, are the electrical resistivity,
thermoelectric and thermal conductivity tensors, respectively.  From
these quantities, we are able to compare the heat flux and electric
field issued from the Fokker-Planck solver to those calculated
analytically in \cite{epperlein}, in the local regime.\\
The classical derivation procedure to obtain the transport
coefficients involves the linearization of the Fokker-Planck-Landau
equation, assuming the plasma to be close to the thermal
equilibrium. The distribution function is approximated using a
truncated Cartesian tensor expansion $\displaystyle f(t,{\bf x},{\bf
v})=f^{(0)}(\|{\bf v}\|^2)+\frac{\bf v}{\|{\bf v}\|^2} \cdot {\bf
f^{(1)}}(t,{\bf x},{\bf v})$. Following \cite{epperlein}, ${\bf \Pi}$
and $m_en_e<{\bf v}>^2$ are neglected. Then considering appropriate
velocity moments of ${\bf f^{(1)}}$, electric fields and heat fluxes
are expressed as a function of thermodynamical variables. The
coefficients of proportionality, in the obtained relations, are
defined as the transport coefficients.\\ Several notations can be
used, depending on the chosen thermodynamical variables. Adopting the
Braginskii notations, we obtain

\begin{eqnarray} 
 \label{braginskii2}
\left\{
\begin{array}{l}
\displaystyle {\bf R} = {\bf \nabla} p + en_e{\bf E} - {\bf j} \times {\bf B} =
\frac{{\bf \alpha_{ep}} \cdot {\bf j}}{n_ee} - {\bf \beta_{ep} \cdot \nabla}T_e, \\ \, \\
 \displaystyle {\bf q}= -\frac{5}{2}\frac{\bf
j}{e}T_e - {\bf \kappa_{ep}}\cdot{\bf \nabla}T_e- {\bf
\beta_{ep}}\cdot {\bf j} \frac{\bf T_e}{e} .
\end{array}\right.
\end{eqnarray}
We want to compare of the results of the solver with the analytical
electric fields and heat fluxes in the local regime. For that purpose,
we use the values of coefficients, for $Z=1$, that are tabulated in
\cite{epperlein}. As for the components of these tensors, we make use
of the standard notations $||$, $\perp$, and $\wedge$. Directions
denoted with $||$ and $\perp$ are respectively parallel and
perpendicular to the magnetic field. Consequently, the parallel and
perpendicular components of a vector ${\bf u}$ are respectively $
u_{||}={\bf b} ({\bf u} \cdot {\bf b})$ and $ u_{\perp} = {\bf b}
\times ({\bf b} \times {\bf u}) $, where ${\bf b}$ is the unit vector
in the direction of the magnetic field. The direction defined by the
third direction in a direct orthogonal frame is denoted by
$\wedge$. In the system \eqref{braginskii2}, the relation between any
transport coefficient tensor ${\bf \varphi}$ and vector ${\bf u}$ is
defined by
\begin{eqnarray} 
 \label{relation1}
 {\bf \varphi}\cdot {\bf u} = \varphi_{||}{\bf b}({\bf b\cdot u}) +
\varphi_{\perp}{\bf b} \times ( {\bf u} \times {\bf b}) \pm
\varphi_{\wedge}{\bf b} \times {\bf u} \ ,
\end{eqnarray}
where the negative sign applies only in the case ${\bf \varphi}={\bf
\alpha_{ep}}$. These coefficients can be expressed in dimensionless
form
\begin{eqnarray} 
 \label{dimbraginskii2}
\left\{
\begin{array}{l}
 \displaystyle {\bf \alpha^c_{ep}} ={\bf
\alpha_{ep}}\frac{\tau}{m_en_e}, \\ \, \\ \displaystyle {\bf
\beta^c_{ep}} = {\bf \beta_{ep}}, \\ \, \\ \displaystyle {\bf
\kappa^c_{ep}} = {\bf \kappa_{ep} }\frac{m_e}{n_e \tau T_e}.
\end{array}\right.
\end{eqnarray}
The dimensionless transport coefficients ${\bf \alpha_{ep}^c}$, ${\bf
\beta_{ep}^c}$, ${\bf \kappa_{ep}^c}$ are functions of $Z$ and the
Hall parameter $ \chi= \omega_c \tau$ only.\\ The heat flux and the
electric field in (\ref{braginskii2}) can then be rewritten in terms
of dimensionless quantities, for the particular 1D geometry of our
temperature gradient configuration. In that case, the normalizations
using a collision frequency (\ref{eq_def_n2})-(\ref{scalingb1}) are
used.
\begin{eqnarray} 
 \label{scal1bis}
\left\{
\begin{array}{l}
 \displaystyle q_1 = -\frac{5}{2}T_en_e^{-1}j_1 -\chi T_e B_3^{-1}
\nabla_{x_1} T_e {\bf \kappa_{ep, \perp}^c} - T_e \left ( {\bf
\beta_{ep, \perp}^c} j_1 - {\bf \beta_{ep, \wedge}^c} j_2 \right ), \\
\, \\ \displaystyle q_2 = -\frac{5}{2}T_en_e^{-1}j_2 -\chi T_e
B_3^{-1} \nabla_{x_1} T_e {\bf \kappa_{ep, \wedge}^c} - T_e \left (
{\bf \beta_{ep, \perp}^c} j_2 + {\bf \beta_{ep, \wedge}^c} j_1 \right
), \\ \, \\ \displaystyle E_1 = {n_e}^{-1} j_2B_3 -
n_e^{-1}\nabla_{x_1} p - \nabla_{x_1} T_e {\bf \beta_{ep, \perp}^c}
+{n_e}^{-1}B_3{\chi}^{-1}({\bf \alpha_{ep, \perp}^c} j_1 + {\bf
\alpha_{ep, \wedge}^c} j_2), \\ \, \\ \displaystyle E_2 =
-{n_e}^{-1}j_1B_3 - \nabla_{x_1} T_e {\bf \beta_{ep, \wedge}^c}
+{n_e}^{-1}B_3{\chi}^{-1}( {\bf \alpha_{ep, \perp}^c} j_2 - {\bf
\alpha_{ep, \wedge}^c} j_1).
\end{array}\right.
\end{eqnarray}
The Hall parameter $ \chi$ is expressed in terms of the dimensionless
quantities $B_3$ and $T_e$:
\begin{eqnarray} 
 \label{hallscal1bis}
\chi = \frac{3\sqrt{\pi}}{2\sqrt{2}}\frac{B_3T_e^{3/2}}{ Z}.
\end{eqnarray}

\subsubsection{Electron transport in the local regime}

In order to validate the numerical scheme in the local regime, we
compare the heat flux ${\bf Q_{FP}}$ and electric field ${\bf E_{FP}}$
computed from the numerical solution, with those analytically (denoted
by ${\bf Q_{BR}}$ and ${\bf E_{BR}}$) computed from the system
(\ref{scal1bis}). The transport coefficients ${\bf \alpha_{ep}}$,
${\bf \beta_{ep}}$, ${\bf \kappa_{ep}}$ have been tabulated in
\cite{epperlein}.\\

In this test case the source term can be considered stiff; the
discretization of the collision operator is then of crucial importance
and its accuracy can be tested. Moreover we provide, in this local
regime, with validation results for a wide range of Hall parameters
corresponding to ICF applications.\\
The initial temperature gradient $T_e(x_1)$ has the form of a step

\begin{eqnarray} 
 \label{Te}
T_e(x_1) = \left\{
\begin{array}{l}
 \displaystyle T_e^R(x_1) \qquad \mbox{if} \quad x_1>x_1^m \ , \\ \, \\ \displaystyle
 T_e^L(x_1) \qquad \mbox{else} \ ,
\end{array}\right.
\end{eqnarray}
 where $T_e^R$ and $T_e^L$ are third order polynomials in
$x_1-x_{1}^{m}$, $x_1$ standing for the space coordinate and
$x_{1}^{m}$ for the mid-point of the 1D domain. The coefficients of
these polynomials are chosen such as they verify the following
conditions at $x_1^m$
\begin{eqnarray} 
 \label{poly}
\left\{
\begin{array}{l}
 \displaystyle \frac{\partial T_e^L}{\partial x_1}(x_{1}^{m})=
\frac{\partial T_e^R}{\partial x_1}(x_{1}^{m})=
\frac{T_R-T_L}{(x_{1}^{R}-x_{1}^{L})/\lambda}, \\ \, \\ \displaystyle
T_e^L(x_{1}^{m})= T_e^R(x_{1}^{m})=\frac{T_R+T_L}{2},
\end{array}\right.
\end{eqnarray}
and at the boundaries
\begin{eqnarray} 
 \label{poly2}
\left\{
\begin{array}{l}
\displaystyle T_e^L(x_{1}^{L})=T_L, \\ \,\\ \displaystyle
T_e^R(x_{1}^{R})=T_R, \\ \,\\ \displaystyle \frac{\partial
T_e^L}{\partial x_1}(x_{1}^{L})= \frac{\partial T_e^R}{\partial
x_1}(x_{1}^{R})=0,
\end{array}\right.
\end{eqnarray}
where $T_L$ (resp. $T_R$) is the initial temperature of the leftmost
(resp. rightmost) point $x_{1}^L$ (resp. $x_{1}^R$) of the
domain. $\lambda$ is a parameter that determines the initial stiffness
of the temperature gradient.\\ \\

The simulations were performed with the following parameters: the
uniform magnetic field $B_3(t=0,x_1)=0.001 , 0.01 , 0.1 , 1$, the size
of the dimensionless domain $L=x_{1}^R-x_{1}^L=5400$, $2 \times
v_{max}=12$, the ion charge $Z=1$, the frequency ratio
$\nu_{e,i}/\omega_{pe}=0.01$, the electron thermal velocity such as
$v_{th}/c=0.05$. The initial electric field is zero over the domain:
$E_1(t=0,x_1)=E_2(t=0,x_1)=0$. The initial distribution function is a
Maxwellian depending on the local temperature, the density being
constant over the domain. The initial temperature profile is chosen
such as $T_L=0.8$, $T_R=1.2$ and $\lambda=10$. This set of parameters
enable us to consider the local regime, close to the hydrodynamics
(the Knudsen number is about $1/500$). The dimensionless time step and
meshes size are $\Delta t=1/500$, $\Delta x_1 = L/126$, $\Delta v =
2v_{max}/32$ respectively. The grid has $126$ points in space and
$32^3$ points in velocity; $42$ processors were used for each
simulation (CEA-CCRT-platine facility). Domain decomposition on the
space domain allows each processor to deal only with $3$ points in
space. The fourth order scheme on the transport part has been used.\\
\\

\begin{figure}[htb]
       \centering
 \includegraphics[angle=0,height=7.5cm,width=7.5cm,keepaspectratio=true]{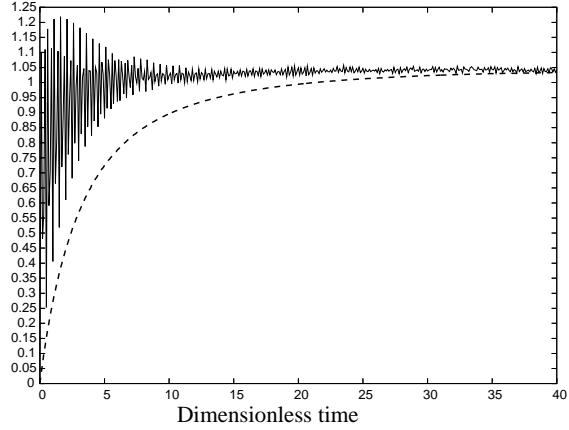}
\caption[]{{\it Longitudinal (along the temperature gradient) ratios
$\frac{ \max_{x_1}(Q_{FP}) }{ \max_{x_1}(Q_{BR}) }$ (dashed curve) and
$\frac{ \max_{x_1}(E_{FP}) }{ \max_{x_1}(E_{BR}) }$ (oscillating
curve) are shown against the dimensionless time. The dimensionless
magnetic field is $B_3=0.001$. Asymptotic behaviour, where the flux is
well established, shows good agreement (less than $5\%$ error) with
analytical solution (Braginskii formalism), denoted by subscript
$\small{BR}$. \label{fig2}}}
\end{figure}

\begin{figure}[htbp]
\begin{tabular}{cc}
\includegraphics[width=4.5cm,height=7.25cm,angle=-90]{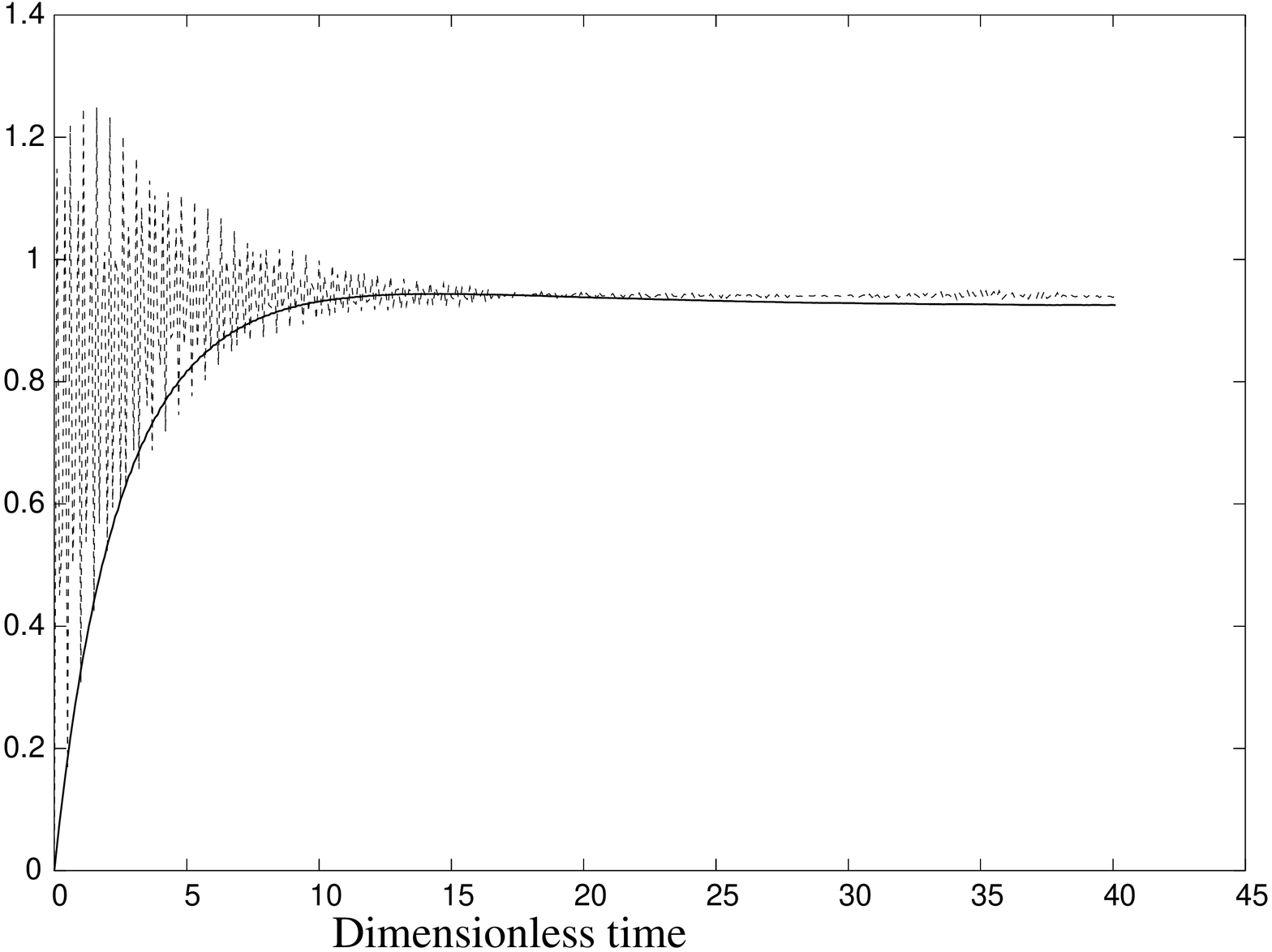}
&
\includegraphics[width=4.5cm,height=7.25cm,angle=-90]{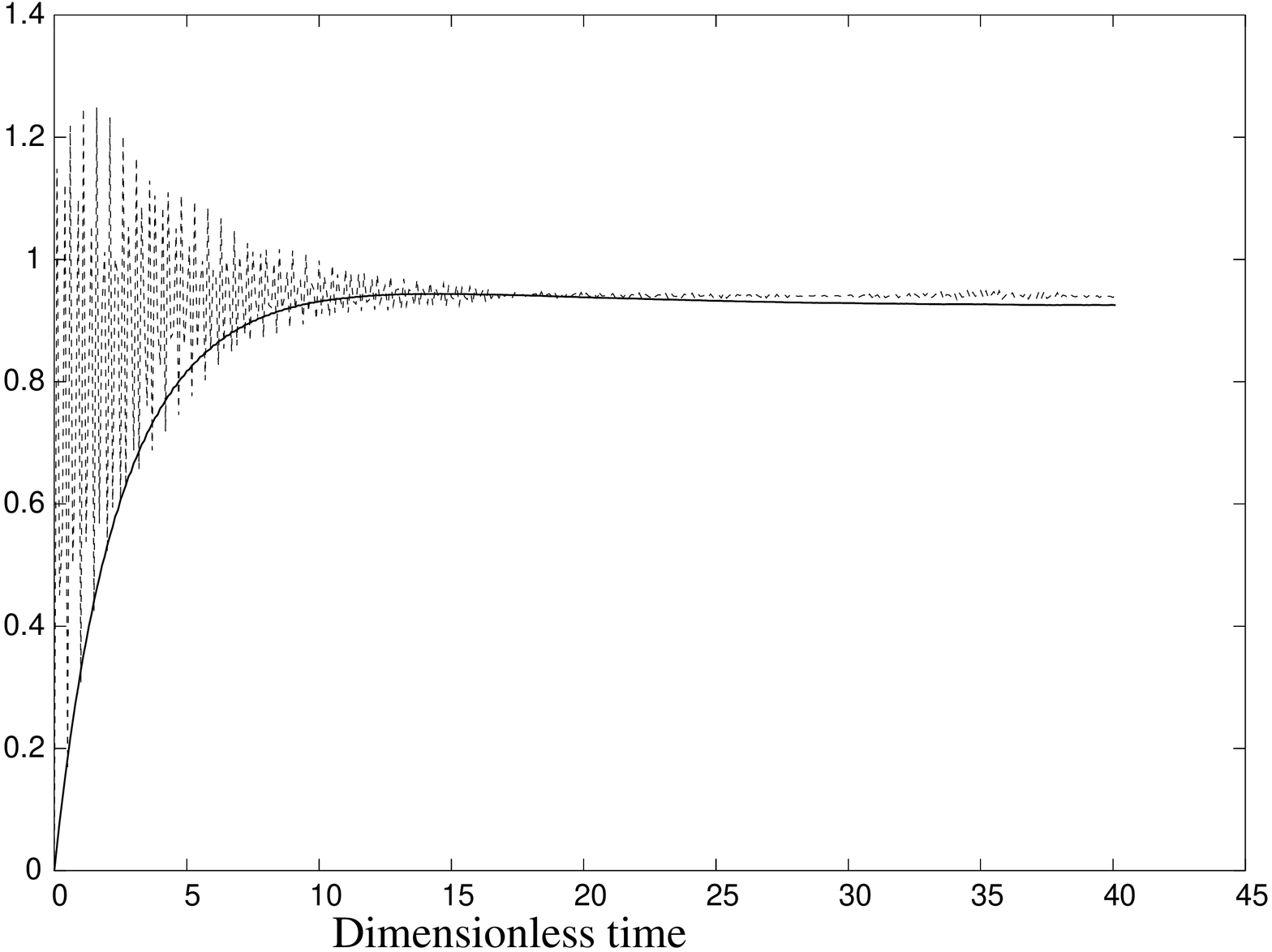}
\\ (a) Longitudinal &(b) Transverse
\end{tabular}
 \caption{{\it Ratios $\frac{ \max_{x_1}(Q_{FP}) }{ \max_{x_1}(Q_{BR})
}$ (curve in bold) and $\frac{ \max_{x_1}(E_{FP}) }{
\max_{x_1}(E_{BR}) }$ (dashed curve) are shown against the
dimensionless time. Longitudinal quantities (along the temperature
gradient) are shown in (a), with about $10\%$ accuracy in the
asymptotics. Transverse quantities are shown in (b), with about $20\%$
accuracy in the asymptotics. The dimensionless magnetic field is
$B_3=0.1$. \label{fig2bis}}}
\end{figure}

\begin{figure}[htbp]
\begin{tabular}{cc}
\includegraphics[width=7.25cm,height=4.5cm,angle=0]{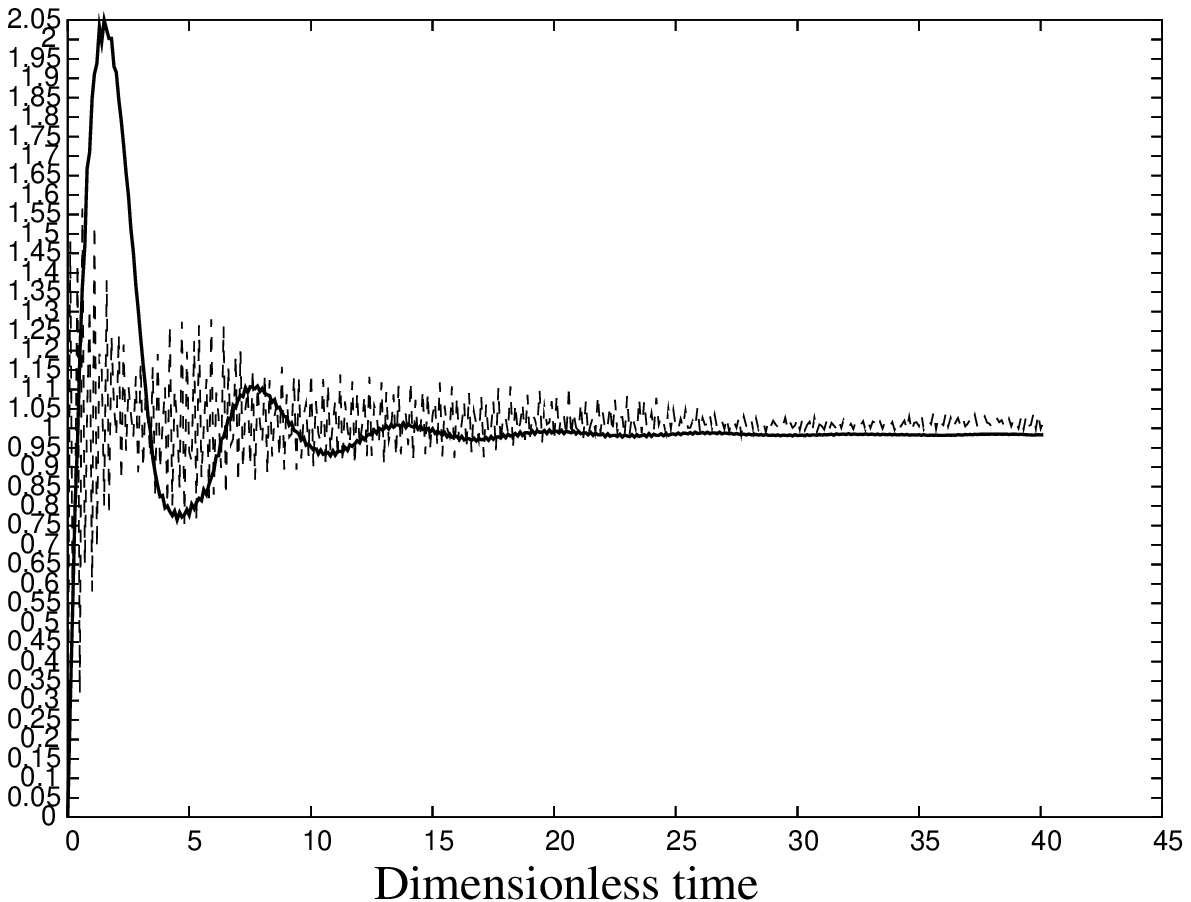}
&
\includegraphics[width=7.25cm,height=4.5cm,angle=0]{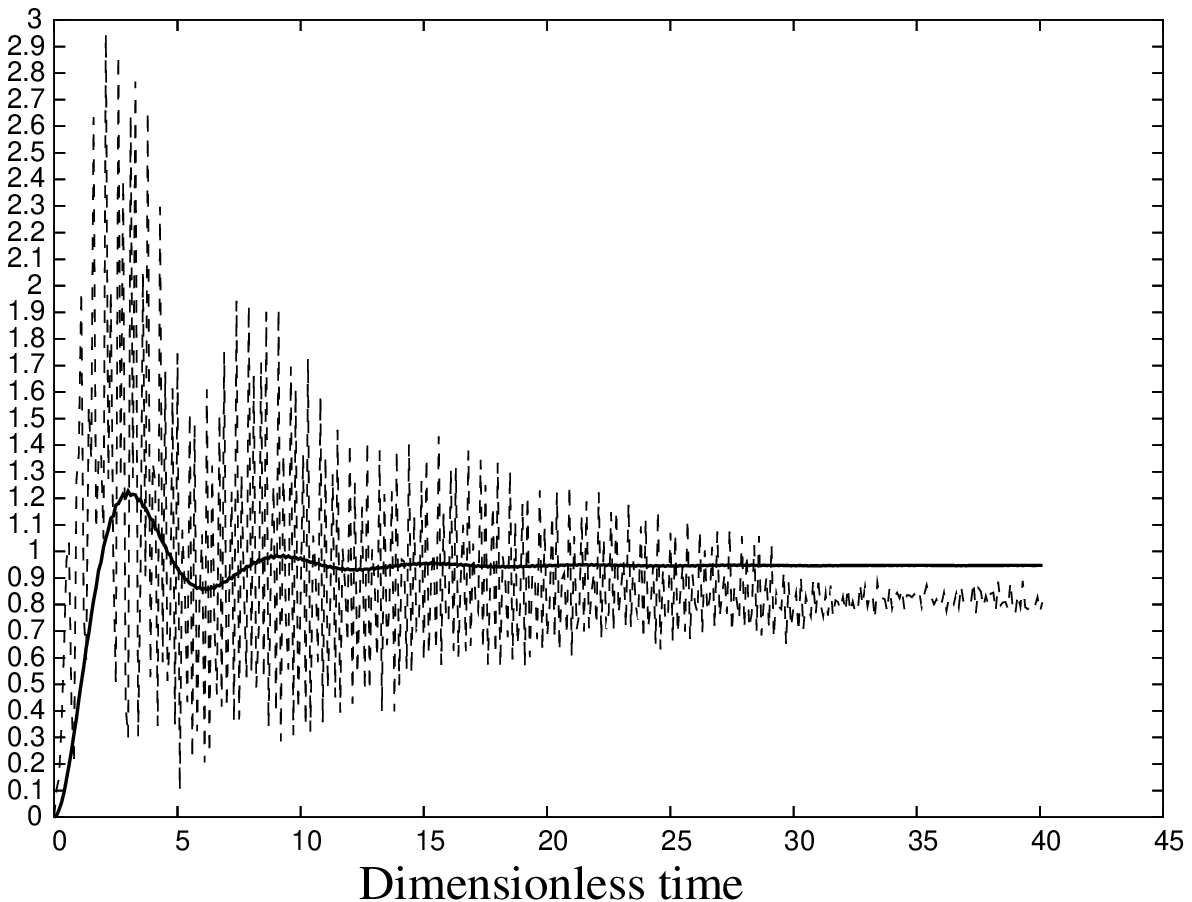}
\\
(a) Longitudinal &(b) Transverse
\end{tabular}
 \caption{{\it Ratios $\frac{ \max_{x_1}(Q_{FP}) }{\max_{x_1}(Q_{BR})
}$ (curve in bold) and $\frac{ \max_{x_1}(E_{FP}) }{
\max_{x_1}(E_{BR}) }$ (dashed curve) are shown against the
dimensionless time. Longitudinal quantities (along the temperature
gradient) are shown in (a), with about $5\%$ accuracy in the
asymptotics. Transverse quantities are shown in (b), with about $20\%$
accuracy in the asymptotics. The dimensionless magnetic field is
$B_3=1$. \label{fig3}}}
\end{figure}

  Results are presented in Figures~\ref{fig2}-\ref{fig3}. The typical
run time is 24 hours for 40 collision times, with that set of
parameters. The maximum difference between the numerical and the
analytical solution are less than $ 10\% $ for longitudinal
macroscopic quantities (heat flux and electric field); $ 20\% $ for
transverse ones. Transverse quantities have only been considered for
simulations presented in Figures~\ref{fig2bis} and \ref{fig3} where
the magnetic field was strong enough so that
\begin{itemize}
\item The establishment of transverse heat flux can be acheived during the simulation time.   
\item Transverse quantities cannot be considered negligible compared to longitudinal ones. 
\end{itemize}

These conditions where fulfilled for $B_3=0.1,1$.\\ In
Figures~\ref{fig2}-\ref{fig3}, only results for simulations with
$B_3=0.001$, $B_3=0.1$, $B_3=1$ are shown, respectively. The
simulation with $B_3=0.01$ proved to show no significance differences
with those with $B_3=0.001$.

Results shown Figures~\ref{fig2}-\ref{fig3} are revealing an important
transient phase before the establishment of a stationary regime. The
oscillations are enforced by the magnetic field, Figure~\ref{fig3}. The oscillating electric fields are the consequence of
the plasma waves excited by our initial conditions; they are damped in
a few electron-ion collision times. These plasma oscillations are
smoothed out by the large time steps we used in simulations, allowed
by the implicit treatment of the Maxwell equations. However this has a
little importance on the asymptotic values and a little importance for
accuracy. With a larger magnetic field Figure~\ref{fig3}, we observe
frequency modulations at $\omega_c=\nu_{e,i}$ (corresponding to
$B_3=1$), both on electric fields and heat fluxes.\\

 In order to investigate Larmor radius effects for simulations
presented in Figures~\ref{fig2bis} and \ref{fig3}, we refined the
space grid below the dimensionless Larmor radius
$r_L=B_3^{-1}$. Therefore, simulation presented in
Figure~\ref{fig2bis} has been done again with the same parameters on
the same time period: we have refined the grid to $1260$ points in
space (420 processors). In the same manner, the simulation presented
in Figure~\ref{fig3} has been done again with $6300$ grid points in
space (2100 processors) and $\Delta t =1/1000$ (C.F.L. condition),
during the same time period. The results prove to be similar to those
with coarse space grids, both for macroscopic quantities and
distribution functions. We thus show no dependence on the Larmor
radius. Here we remark that the cyclotron frequency is always
resolved. The time steps are constrained, for most of the cases we
treat, by the C.F.L. on collision operators.

\subsubsection{Electron transport in the nonlocal regime}

The departure of transport coefficients from their local values is of
interest here. We restrict ourselves to cases where the magnetic
field is zero. Then it is possible to obtain directly the ratio of
effective thermal conductivity to the Spitzer-H\"{a}rm conductivity
$\kappa / \kappa_{SH}$ by the relation:
\begin{eqnarray} 
 \label{effectconductivity}
 \frac{\kappa}{\kappa_{SH}} = \frac{q_1}{q_{SH}}.
\end{eqnarray}
The Spitzer-H\"{a}rm regime refers to a local regime with no magnetic
field.  In (\ref{effectconductivity}), $q_1$ is calculated from the
numerical solution and $q_{SH}$ from (\ref{scal1bis}) in the
Spitzer-H\"{a}rm limit.\\ Transport coefficients are extracted from
the domain where the flux and temperature gradient are maximum.\\ The
wavelength of the temperature perturbation $k\lambda_{e,i}$ in the
Fourier space is computed from the gradient temperature profile. This
enables to obtain a range (due to an uncertainty) for $k\lambda_{e,i}$
corresponding to this temperature gradient.\\ The results will be
compared with the analytical formula from \cite{epperlein2}

\begin{eqnarray} 
 \label{af1}
&& \frac{\kappa}{\kappa_{SH}} = \frac{1}{1+(30k\lambda_{e,i}
\beta)^{4/3}}, \\ && \beta = \left (
\frac{3\pi}{128}\frac{3.2(0.24+Z)}{(1+0.24Z)}\right)^{1/2}
\frac{Z^{1/2}}{2}.
\end{eqnarray}

The comparison between the numerical results and the analytical
solution are in good agreement. The three runs have been performed
with the same precision for the temperature gradient.\\

\begin{tabular}{|l|c|c|c|}
\hline \qquad \mbox{\bf Parameters} & RUN1 & RUN2 & RUN3 \\ \hline
Size of the domain & 5400 & 540 & 540 \\ \hline Stiffness parameter
$\lambda$ & 10 & 10 & 100 \\ \hline Number of points along the
Gradient & 126 & 126 & 1260 \\ \hline Number of processors & 42 & 42 &
420 \\\hline \qquad \mbox{\bf Results} & RUN1 & RUN2 & RUN3 \\ \hline
\qquad $k\lambda_{e,i}$ & $ 10^{-3}$ & $0.05 \pm 0.03$ & $0.2 \pm 0.1$
\\ \hline Analytical $\kappa / \kappa_{SH} $ & $0.998$ & $[0.93-0.67]$
& $[0.60-0.26]$ \\ \hline Numerical $\kappa / \kappa_{SH} $ & 1.03 &
0.675 & 0.395\\ \hline
\end{tabular}

\subsection{2D nonlocal magnetic field generation}

We present here results on the nonlocal magnetic field generation
during the relaxation of cylindrical laser hot spots, having a
periodic repartition, and for a region of constant density. This
stands as a first step to prove the $2D$ capabilities of the
solver. The 2D extension of the presented numerical schemes is
straightforward on a grid.\\

We consider a planar geometry with periodic boundary conditions. For
this application, the normalizations using collision frequency
(\ref{eq_def_n2})-(\ref{scalingb1}) are used.\\ The initial
dimensionless temperature profile is $T_e({\bf x},t=0)=1+0.12 \exp
\left ( -\frac{{\bf x}^2}{R^2}\right ) $ , with $R=5.6$. We used the
following parameters for the simulation: the frequency ratio is
${\nu_{e,i}}/{\omega_{pe}} = 0.003 $, the ion charge $Z$ is assumed to
be high, so that we do not consider the electron-electron collision
operator; here the relaxation only acts with electron-ion collisions
on the anisotropic part of the electronic distribution function.  The
electron thermal velocity is such as $v_{th}/c=0.05$. These parameters
are close to those used in \cite{senecha}. The size of the simulation
domain is $L=70$ for one space direction, $2 \times v_{max}=12$ for
one velocity direction. Initial electric and magnetic fields are zero
over the domain. The initial distribution function is a Maxwellian
depending on the local temperature, the density being constant over
the domain. The dimensionless time step and meshes size are $\Delta
t=1/500$, $\Delta x=\Delta y= L/100$, $\Delta v= 2v_{max}/32$,
respectively. The grid has $100^2$ points in space and $32^3$ points
in velocity. $625$ processors are used for this simulation.

\begin{figure}[htbp]
\begin{tabular}{cc}
\includegraphics[width=9.5cm,height=8.25cm,angle=-90]{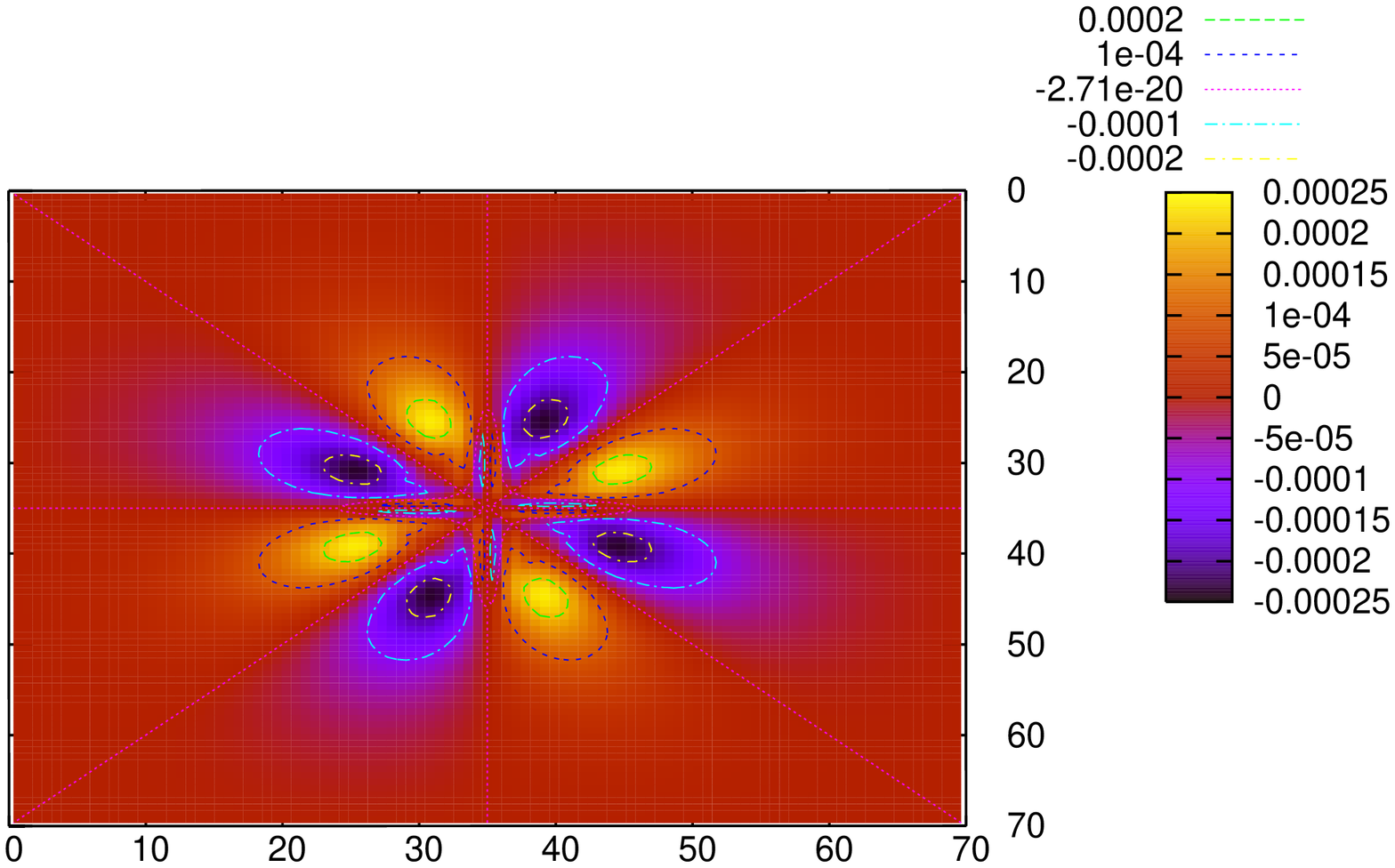}
&
\includegraphics[width=9.5cm,height=8.25cm,angle=-90]{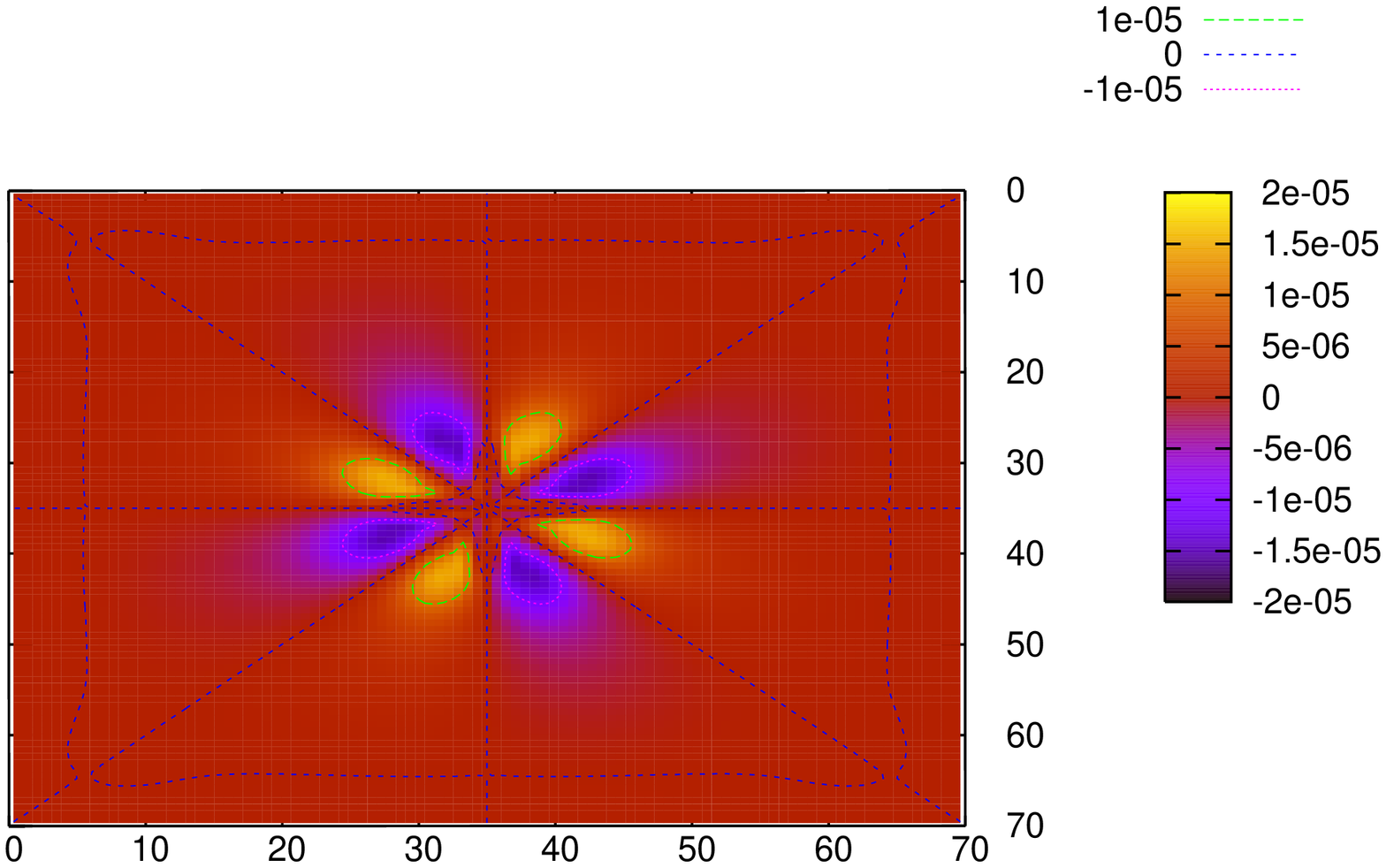}
\\
(a) Magnetic field &(b) Cross gradients of high order moments
\end{tabular}
 \caption{{\it Dimensionless magnetic field and cross gradients of high
 order moments (third and fifth) at $t\nu_{e,i}=8$. \label{final}}}
\end{figure}
The mecanism under consideration here (the magnetic field generation
Figure~\ref{final}), is expained in \cite{Kingham}, as the results of
non parallel gradients of the third and fifth moments of the
electronic distribution function. We show the magnetic field in
Figure~\ref{final},(a) and the cross gradients $\displaystyle
\nabla_{\bf x} \left ( \int_{\RR^3} f_e \| {\bf v} \|^3 d{\bf v}
\right ) \times \nabla_{\bf x} \left ( \int_{\RR^3} f_e \| {\bf v}
\|^5 d{\bf v} \right ) $ in Figure~\ref{final},(b).\\ This mecanism
is not due to the magnetic field generation from a $\nabla n_e \times
\nabla T_e$ structure, since the density $n_e$ remains constant over
the domain.\\ This structure with eight lobes is the result of the
collision operators (of diffusion type) that make a particular speckle
interact, after a rapid transient phase, with the other surrounding
(similar) speckles. We note that an important parameter to anayse
further such interactions should be the size of the speckle over the
distance between speckles.
\section{Conclusions}

In the present paper, we have developed high order numerical methods
dedicated to plasma simulation at a microscopic level.

A fourth order scheme issued from VFRoe schemes has been introduced in
our kinetic context. It brings accuracy improvement on the velocity
transport term. The second order scheme remains interesting for the
linear spatial transport term (which is faced to less robustness and
accuracy constraints) in a 2D, distributed memory context without
overlapping between processors (each processor communicating with its
neighbours only). It involves indeed a reduced stencil allowing for a
lower minimum number of spatial grid points per processor.\\ The
Maxwell equations have been discretized with a second order, implicit
scheme allowing large time steps. We did not find any dependance on
the Larmor radius and show that resolving the cyclotron frequency is
sufficient. The couplings between the equations of the system have
introduced a number of constraints (robustness, accuracy, symmetry)
both on the transport scheme and the collision operators. Some
numerical and physical test cases have validated our approach in
different regimes of interest for ICF applications, and showed that it
is computationally affordable. We also proposed a validation strategy
in the linear regime based on \cite{sc}, using Green kernels. \\
Various fundamental studies can be planned on the basis of the actual
version of the solver. Collisional Weibel instability \cite{sangam},
forward and backward collisional Stimulated Brillouin Scattering,
studies on the nonlocal interaction beween speckles for plasma-induced
smoothing of laser beams issues \cite{feugeas}, for instance. Also
several axis of development are under consideration to bring more
physics to the model: the ion motion, the extension to regimes
relevant to higher laser intensities (relativistic regime and large
angle collision terms of Boltzmann type).

\paragraph{\bf Acknowledgments:}
The authors are thankful to the Commissariat à l'Energie Atomique for
the access to the CEA-CCRT-platine computing facillities. One of the
author, Francis Filbet, would like to express his gratitude to the ANR
JCJC-0136 MNEC (M\'{e}thode Num\'{e}rique pour les Equations
Cin\'{e}tiques) funding.

%\end{document}

\appendix
\section{Electrostatic case in the linear regime}
\label{annexeA}

The relativistic $1D_x \times 3D_{\bf v}$Vlasov-Poisson system extracted
from the equations \eqref{MFP1}-\eqref{MFP5} reads
\begin{eqnarray} 
  \label{vla_rel}
  && \frac{\partial f}{\partial t } + v_1 \frac{\partial f}{\partial x_1}  + \frac{q_e}{m_e} {E_1} \frac{\partial f}{\partial v_1 } =  0 \ ,
  \\
\label{poisson}
  && \frac{\partial E_1}{\partial x_1} = -\frac{q_e}{\epsilon_0} \left ( n_0 -
  \int_{\RR^3} f(t,x_1,{\bf v}){\bf dv} \right ) \ .
\end{eqnarray}
The distribution function $f$ is assumed to be a perturbation around
an equilibrium state $f^{(0)}(\|{\bf v}\|)$, $E^{(0)}_1=0$, $n_{0} = \displaystyle
\int_{\RR^3} f^{(0)}(\|{\bf v}\|) {\bf d{v}}$. The system
(\ref{vla_rel}),(\ref{poisson}) is linearized around this equilibrium
state
\begin{eqnarray} 
  \label{distrib}
&& f(t,x_1,{\bf v})=f^{(0)}(\|{\bf v}\|)+f^{(1)}(t,x_1,{\bf v}) \ ,
\\
&& E^{(1)}_1(t,x_1)=E^{(0)}_1+E^{(1)}_1(t,x_1) \ ,
\end{eqnarray}
under the hypothesis:
\begin{eqnarray} 
 \label{inf0}
\|f^{(1)}\| \ll \|f^{(0)}\| \ , 
\\
\label{inf1}
 \|E^{(1)}_1\| \ll 1 \ .
\end{eqnarray}
The Vlasov-Poisson can then be set under the following form (transport
equation along the space directions supplemented by a source term
along the $v_1$ direction), after linearization
\begin{eqnarray} 
  \label{linearise}
 \left\{ \begin{array}{ll} \displaystyle \frac{\partial
      f^{(1)}}{\partial t} + v_1 \frac{\partial f^{(1)}}{\partial x_1}
      = -\frac{q_e}{m_e} E^{(1)}_1 \frac{\partial f^{(0)}}{\partial
      v_1} \ , \\ \, \\ \displaystyle \frac{\partial
      E^{(1)}_1}{\partial x_1} = \frac{q_e}{\epsilon_0} \int_{\RR^3}
      f^{(1)}(t,x_1,{\bf v}){\bf dv} \ .
 \end{array}\right.
\end{eqnarray}
If $f^{(1)}$ and $E^{(1)}_1$ are periodic and integrable, then their
respective normalized Fourier coefficient are well-defined. A Fourier
series expansion gives $\forall t >0$ 
\begin{eqnarray}
\label{perturbation3}
 \left\{ \begin{array}{ll} \displaystyle f^{(1)}(t,x_1,{\bf v})=
 \hat{f}^{(1)}(t,k_1,{\bf v})cos(k_1x_1) \ , \\ \, \\ \displaystyle
 \hat{f}^{(1)}(t,k_1,{\bf v})=\frac{1}{L} \int_{0}^L
 f^{(1)}(t,x_1,{\bf v})e^{-ik_1x_1}dx_1 \ ,
 \end{array}\right.
\end{eqnarray}
Where $L$ is the size of the domain. The same reconstruction using
Fourier series is used for $E^{(1)}_{1}$.\\ These coefficients verify
the following equations,obtained by Fourier transformation performed
on the equations of the system \eqref{linearise}, for all real $k_1$
\begin{eqnarray} 
  \label{sysfou}
&& \frac{\partial \hat{f}^{(1)}}{\partial t}+ ik_1v_1 \hat{f}^{(1)}=
-\frac{q_e}{m_e} \hat{E}^{(1)}_1 \frac{\partial f^{(0)}}{\partial v_1} \ ,
\\
\label{sysfou2}
&& ik_1 \hat{E}_1=\frac{q_e}{\epsilon_0} \hat{n}_1 \ .
\end{eqnarray}
Then introducing the notation $\hat{f}^{(1)}(t=0,k_1, {\bf v}
)=\hat{f}^{(10)}(k_1,{\bf v})$, the equation (\ref{sysfou}) can be
written in the integral form
\begin{eqnarray} 
  \label{sol1}
&& \hat{f}^{(1)}(t,k_1,{\bf v})= \hat{f}^{(10)}(k_1,{\bf v})
      e^{-ik_1v_1t} - \frac{q_e}{m_e} \int_{0}^{t} \hat{E}^{(1)}_1(t',k_1)
      \frac{\partial f^{(0)}}{\partial v_1} e^{-ik_1v_1(t-t')}dt' \ .
\end{eqnarray}
Integrating the equation \eqref{sol1} over ${\bf v}$ and injecting in
it the relation \eqref{sysfou2}, one obtains the following integral
equation for the density
\begin{eqnarray}
  \label{sol3}
\hat{n}^{(1)}(t,k_1)= M(t,k_1)+ \int^{t}_{0} K(t-t',k_1)\hat{n}^{(1)}(t',k_1)dt' \ ,
 \end{eqnarray}
where
\begin{eqnarray} 
  \label{sol4}
&& K(t,k_1)= \frac{iq_e^2}{k_1m_e\epsilon_0} \int_{\RR^3}
    \frac{\partial f^{(0)}}{\partial v_1}e^{-ik_1v_1t}{\bf d{v}} \ ,
\\
&& M(t,k_1)= \int_{\RR^3} \hat{f}^{(10)}(k_1,{\bf
    v})e^{-ik_1v_1t}{\bf d{v}} \ .
 \end{eqnarray}
These kernels can be computed with the desired accuracy, following
\cite{sc}.
The numerical resolution of \eqref{sol3} finally reduces to the
inversion of a triangular linear system.\\  
Macroscopic quantities such as the density or the heat flux can then
be reconstructed using these latter equations.

\section{Initialisation for the generation of a single X-mode plasma wave}
\label{annexeC}

This test case stands as a validation for the couplings of Vlasov and
Maxwell equations. We determine initial conditions that trigger a
plasma wave at a given wavelength. To do so, Vlasov-Maxwell equations
are linearized, setting $f=f^{(0)}+\tilde{f}$, $E=\tilde{E}$,
$B=B^{(0)}+\tilde{B}$ around the equilibrium state $f= f^{(0)}$, $E =
0$, $B = B^{(0)}$. In this appendix, we use the normalization
\eqref{scaling}-\eqref{scaling1}. We assume periodic boundary
conditions. The fluctuations of the total pressure tensor are
neglected with respect to those of the magnetic field.\\
Using the conservation law $\displaystyle \frac{\partial n}{\partial t} + \frac{\partial j_1}{\partial {x_1}} =0$, the
former hypothesis lead us to solve the system of six equations with
six unknown $\tilde{j}_1$, $\tilde{j}_2$, $\tilde{E}_1$,
$\tilde{E}_2$, $\tilde{B}_3$ and $\tilde{n}$

\begin{equation}
  \label{sys_dispersif}
  \left\{ \begin{array}{ll}
      \displaystyle
      \frac{\partial \tilde{j}_1}{\partial t} + \tilde{E}_1 + B^{(0)} \tilde{j}_2 =0 \ , \\
      \displaystyle
      \frac{\partial \tilde{j}_2}{\partial t} + \tilde{E}_2 - B^{(0)} \tilde{j}_1 =0 \ , \\
      \displaystyle
      \frac{\partial \tilde{n}}{\partial t} +\partial_{x_1} \tilde{j}_1 =0 \ , \\
      \displaystyle
      \frac{\partial \tilde{E}_1}{\partial x_1 } = -\tilde{n} \ , \\
      \displaystyle
      \frac{\partial\tilde{E}_2}{\partial t} = -
      \frac{1}{\beta^2} \frac{\partial \tilde{B}_3}{\partial x_1} +
      \tilde{j}_2 \ , \\
      \displaystyle
      \frac{\partial
      \tilde{B}_3}{\partial t} = - \frac{\partial
      \tilde{E}_2}{\partial x_1} \ .
    \end{array}\right.
\end{equation}
Applying time and space Fourier tranform to this system, and
identifying Fourier composants ($\tilde n = \hat n \exp (-i \omega t +
i k_1 x_1)$), the following system is obtain
\begin{displaymath}
  \left\{ \begin{array}{ll}
      \displaystyle  
      -i\omega \hat{j}_{1} + {\hat E}_{1} + B^{(0)} \hat{j}_{2} =0 \ , \\
      \displaystyle  
      -i\omega \hat{j}_{2} + \hat{E}_{2} - B^{(0)} \hat{j}_{1} =0  \ , \\
      \displaystyle  
      -i\omega \hat{n} + ik_1 \hat{j}_{1} =0  \ , \\
      \displaystyle  
      ik_1\hat{E}_{1} = -{\hat n} \ , \\
      \displaystyle  
      -i\omega\hat{E}_{2} = - \frac{1}{\beta^2} ik_1 \hat{B}_{3} +
      \hat{j}_{2}  \ , \\
      \displaystyle  
      -i\omega \hat{B}_{3} = - ik_1 \hat{E}_{2}  \ .
    \end{array}\right.
\end{displaymath}
The dispersion equation of this system reads
\begin{eqnarray}
\label{testtest1}
  N^2 =\frac{k_1^2}{\beta^2 \omega^2} = 1- \frac
  {\omega^2-1}{\omega^2(\omega^2-1-\|B^{(0)}\|^2)} \ .
\end{eqnarray}
In this equation, the plasma frequency is $\omega_{pe}=1 $ and the
cyclotron frequency is $\omega_{c} = q_e\|B^{(0)}\|/m$, that is
$\|B^0\|$ in this dimensionless case. The perturbative term of the
distribution function at initial time can be determined for a
particular solution $\omega$ of this relation dispersion.\\ The
Fourier transform is applied on the linearized Vlasov equation
\begin{equation}
  \label{cinlin}
  (-i \omega+ i k_1 v_1) \hat f - \hat E_1 \frac{\partial f^{(0)}}{
    \partial {v_1} } - \hat E_2 \frac{\partial f^{(0)}}{\partial {v_2}
    } - B^{(0)} v_2 \frac{\partial \hat f}{\partial {v_1} } + B^{(0)}
    v_1 \frac{\partial \hat f}{\partial {v_2}} \,=\,0 \ .
\end{equation}
This equation is expressed in cylindrical coordinates
\begin{displaymath}
  \left\{ \begin{array}{ll}
      \displaystyle
      v_1 = v_{\perp} \cos(\psi) \ , \\
      \displaystyle
      v_2 = v_{\perp} \sin(\psi) \ , \\
      \displaystyle
      v_3 = v_{\parallel}
    \end{array}\right. 
\end{displaymath}
where
\begin{displaymath}
  \left\{ \begin{array}{ll}
      \displaystyle
      v_{\perp}=(\|v_1\|^2 + \|v_2\|^2)^{1/2} \ , \\
      \displaystyle
      \tan (\psi)= \frac{v_2}{v_1} \ .
    \end{array}\right. 
\end{displaymath}
Recalling that:
\begin{gather*}
  \nabla_{\bf v} f = \frac{\partial f}{\partial v_{\perp}} \nabla_{\bf
  v} v_{\perp} +\frac{\partial f}{\partial \psi} \nabla_{\bf v} \psi +
  \frac{\partial f}{\partial v_{\parallel}} \nabla_{\bf v}
  v_{\parallel} \ , \\ \left\{ \begin{array}{ll} \displaystyle
  \frac{\partial v_{\perp}}{\partial v_1}= \cos(\psi) \ , \\
  \displaystyle \frac{\partial v_{\perp}} {\partial v_2}= \sin(\psi) \
  , \\ \displaystyle \frac{\partial \psi}{\partial v_1}=-\frac{1}
  {v_{\perp}} \sin(\psi) \\ \displaystyle \frac{\partial
  \psi}{\partial v_2}=\frac{1} {v_{\perp}} \cos(\psi) \ ,
    \end{array}
  \right. 
\end{gather*}
with $\nabla_{\bf v} v_{\perp} = \vec e_{\perp}$, $\nabla_{\bf v}
v_{\psi} = \vec e_{\psi}$ and $\nabla_{\bf v} v_{\parallel} = \vec
e_{\parallel}$, where $\vec e$ are vectors in the local basis.
Setting $f^{(0)}(\|{\bf v}\|^2) = (2 \pi)^{\frac{3}{2}}
\exp(-\frac{\|{\bf v}\|^2}{2})$, and writing
\begin{gather*}
  ({\bf v} \wedge {\bf B}) . \nabla_{\bf v} \hat f =( \nabla_{\bf v}
  \hat f \wedge {\bf v}) . {\bf B} = - B^{(0)} \frac{\partial
  f}{\partial \psi},\,\, \text {with} \, \, {\bf B}=(0,0,B^{(0)}) \ ,
\end{gather*}
the kinetic equation (\ref{cinlin}) becomes
\begin{equation}
  \label{spher}
  (-i \omega + i k_1 v_{\perp} \cos(\psi)) \hat f +
    B^{(0)}\frac{\partial \hat f}{\partial \psi} + f^{(0)}(\|{\bf
    v}\|^2) v_{\perp} (\hat E_{1} \cos(\psi)+ \hat E_{2}\sin(\psi)) =
    0.
\end{equation}
In order to solve this equation, we decompose the distribution
function as a Fourier serie
\begin{gather*}
  \hat f = \sum_{n=-\infty}^{+\infty} \hat f_n (v_{\perp}) e^{in\psi} .
\end{gather*}
Then from (\ref{spher}),
\begin{gather*}
  \sum_{n=-\infty}^{+\infty} (-i \omega+ i k_1 v_\perp \cos(\psi) + i n
  B^{(0)}) \hat f_n e^{in\psi} = - f^{(0)}(\|{\bf v}\|^2)
  v_{\perp}(\hat E_{1} \cos(\psi)+ \hat E_{2} \sin(\psi)) \ .
\end{gather*}
Multiplying this equation by $e^{i m \psi}$, integrating from $0$ to
$2\pi$, we obtain
\begin{eqnarray}
  \label{syst}
  &&\sum_{n=-\infty}^{+\infty} \int_0^{2\pi} e^{i m \psi} (-i\omega+i
  k_1 v_\perp \cos(\psi) + i n B^{(0)}) \hat f_n e^{i n \psi} d \psi
  \nonumber \\ &&= -f^{(0)}(\|{\bf v}\|^2) v_{\perp} \int_0^{2\pi}
  e^{im\psi} (\hat E_{1} \cos(\psi)+ \hat E_{2} \sin(\psi)) d \psi \ .
\end{eqnarray}
For $m=0$, terms are different from zero only for $n=-1, 0, 1$. From \eqref{syst} comes
\begin{equation}
  \label{syst1}
  k_1v_\perp \hat f_{-1} - 2 \omega \hat f_0 + k_1 v_\perp \hat f_1 =0 \ .
\end{equation}
For $m=-1$,
\begin{equation}
  \label{syst2}
  i k_1 v_\perp \hat f_0 - 2 i (\omega - B^{(0)}) \hat f_1 + i k_1 v_\perp
  \hat f_2 = - f_0(v^2) v_{\perp} (\hat E_{1}- i \hat E_{2}) \ .
\end{equation}
For $m=1$,
\begin{equation}
  \label{syst3}
  i k_1 v_\perp \hat f_{-2} - 2 i (\omega + B^{(0)}) \hat f_{-1} + i k_1
  v_\perp \hat f_0 = - f^{(0)}(\|{\bf v}\|^2) v_{\perp} (\hat E_{1} +
  i \hat E_{2}) \ .
\end{equation}
The case $m=-2$ involves $\hat
f_3$,
\begin{equation}
  \label{syst4}
  i k_1 v_\perp \hat f_1 - 2(\omega - 2 B^{(0)}) \hat f_2 + i k_1 v_\perp
  \hat f_3 = 0 \ .
\end{equation}
In the same manner the case $m=2$ involves $\hat
f_{-3}$,
\begin{equation}
  \label{syst5}
  i k_1 v_\perp \hat f_{-3} - 2(\omega + 2 B^{(0)}) \hat f_{-2} + i k_1
  v_\perp \hat f_{-1} = 0 \ .
\end{equation}
In order to close the system, the components $f_{-3}$ and
$f_3$ are neglected, and we deduce from (\ref{syst1}-\ref{syst5}),
\begin{displaymath}
  \left \{
    \begin{array}{rcl}
      - 2(\omega + 2 B^{(0)}) \hat f_{-2} + i k_1
      v_\perp \hat f_{-1} & = & 0 \ ,
      \\
      i k v_\perp \hat f_{-2} - 2 i (\omega + B^{(0)}) \hat f_{-1} + i k_1
      v_\perp \hat f_0 & = & - f^{(0)}(\|{\bf v}\|^2) v_{\perp} (\hat E_{1} + i \hat
      E_{2}) \ ,
      \\
      kv_\perp \hat f_{-1} - 2 \omega \hat f_0 + k_1 v_\perp \hat f_1
      & = & 0 \ ,
      \\
      i k_1 v_\perp \hat f_0 - 2 i (\omega - B^{(0)}) \hat f_1 + i k_1 v_\perp
      \hat f_2 & = & - f^{(0)}(\|{\bf v}\|^2) v_{\perp} (\hat E_{1}- i \hat E_{2}) \ ,
      \\ 
      i k_1 v_\perp \hat f_1 - 2( \omega - 2 B^{(0)}) \hat f_2 & = & 0 \ .
  \end{array}
  \right.
\end{displaymath}
The solution of linearized Vlasov equation can be calculted
\begin{displaymath}
  \left \{
    \begin{array}{rcl}
      \displaystyle 
      f(t,x,v) &= &f^{(0)}(\|{\bf v}\|^2) + \sum_{n=-\infty}^{+\infty} \hat f_n(v_{\perp})
      e^{-i\omega t + i k_1 x_1 + i n \psi} \ ,
      \\
      \displaystyle 
      E_1(t,x) & = & \hat E_1 e^{-i \omega t + i k_1 x_1} \ ,
      \\
      \displaystyle 
      E_2(t,x) & = & \hat E_2 e^{-i \omega t + i k_1 x_1} \ ,
      \\
      \displaystyle 
      B(t,x)& = & B^{(0)} + \hat B_{3} e^{- i \omega t + i k_1 x_1} \ .
    \end{array}
  \right.
\end{displaymath}
The dispersion relation \eqref{testtest1} provides with a particular
$\omega$. Then we obtain the following results for the construction of the
initial solution,
\begin{displaymath}
  f(0,x,v) = f^{(0)}(\|{\bf v}\|^2) + \sum_{n=-2}^{2} \hat
  f_n(v_{\perp}) e^{i k_1 x + i n \psi}
\end{displaymath}
With the expressions
\begin{displaymath}
  \begin{array}{rcl}
    \displaystyle
    \frac{\hat f_{-2}} {f^{(0)}(\|{\bf v}\|^2) \hat D} & = & 
    i ( -4\,{\omega}^{3}{ \hat E_1}-4\,i{\omega}^{3}{ \hat E_2}+12\,i{
      \omega}^{2}{ B^{(0)}}\,{ \hat E_2}+12\,{\omega}^{2}{ B^{(0)}}\,{ \hat E_1}-8\,{{
        \|B^{(0)}\|}}^{2}\omega\,{ \hat E_1}
    \\   
    \displaystyle
    & + & {k_1}^{2}{{ v_{\perp}}}^{2}\omega\,{ \hat E_1} + 3\,
    i{k_1}^{2}{{ v_{\perp}}}^{2}\omega\,{ \hat E_2}-8\,i{{ \|B^{(0)}\|}}^{2}\omega\,{ 
      \hat E_2}-4\,i{k_1}^{2}{{ v_{\perp}}}^{2}{ B^{(0)}}\,{ \hat E_2} ) {{
        v_{\perp}}}^{2} 
    k_1 , 
    \\
    \displaystyle
    \frac{\hat f_{-1}} {f^{(0)}(\|{\bf v}\|^2) \hat D} & = & 
    2\,i{ v_{\perp}}\, ( { \hat E_1}\,{k_1}^{2}{{ v_{\perp}}}^{2}{\omega}^{2}+4\,i
    { B^{(0)}}\,{\omega}^{3}{ \hat E_2}-16\,{{ \|B^{(0)}\|}}^{3}\omega\,{ \hat E_1}-16\,
    i{{ \|B^{(0)}\|}}^{3}\omega\,{ \hat E_2} 
    \\    
    \displaystyle
     & + & 3\,i{ \hat E_2}\,{k_1}^{2}{{ v_{\perp}}}^{2}{
      \omega}^{2} - 4\,{ \hat E_1}\,{\omega}^{4}-8\,i{ \hat E_2}\,{k_1}^{2}{{ v_{\perp}}}^
    {2}{{ \|B^{(0)}\|}}^{2}+2\,{k_1}^{2}{{ v_{\perp}}}^{2}{ B^{(0)}}\,\omega\,{ \hat E_1}
    \\ 
    \displaystyle
    & + &2 \,i{k_1}^{2}{{ v_{\perp}}}^{2}{ B^{(0)}}\,\omega\,{ \hat E_2}+16\,{ \hat E_1}\,{{
        \|B^{(0)}\|}}^{2}{\omega}^{2} + 16\,i{ \hat E_2}\,{{ \|B^{(0)}\|}}^{2}{\omega}^{2}+4\,
    { B^{(0)}}\,{\omega}^{3}{ \hat E_1}
    \\
    \displaystyle
    & - & 4\,i{ \hat E_2}\,{\omega}^{4} 
    ),
    \\
    \displaystyle
    \frac{\hat f_{0}} {f^{(0)}(\|{\bf v}\|^2) \hat D} & = & 
    2\,i{{ v_{\perp}}}^{2}k_1 ( 16\,{{ \|B^{(0)}\|}}^{2}\omega\,{ \hat E_1}+{k_1}^{2}{
      { v_{\perp}}}^{2}\omega\,{ \hat E_1}-4\,{\omega}^{3}{ \hat E_1}+4\,i{\omega}^{2}
    { B^{(0)}}\,{ \hat E_2} 
    \\ 
    \displaystyle
    & - & 16\,i{{ \|B^{(0)}\|}}^{3}{ \hat E_2} + 2\,i{k_1}^{2}{{ v_{\perp}}}^{
      2}{ B^{(0)}}\,{ \hat E_2} ) ,  
    \\ 
    \displaystyle
    \frac{\hat f_{1}} {f^{(0)}(\|{\bf v}\|^2) \hat D} & = & 
    2\,i ( -2\,{ B^{(0)}}+\omega ) { v_{\perp}}\, ( -4\,i{k_1}^{2}
    {{ v_{\perp}}}^{2}{ B^{(0)}}\,{ \hat E_2}+{k_1}^{2}{{ v_{\perp}}}^{2}\omega\,{ \hat E_1}
    -3\,i{k_1}^{2}{{ v_{\perp}}}^{2}\omega\,{ \hat E_2}
    \\ 
    \displaystyle
    & - & 12\,{\omega}^{2}{ B^{(0)}}\,{
      \hat E_1} + 12\,i{\omega}^{2}{ B^{(0)}}\,{ \hat E_2}-4\,{\omega}^{3}{ \hat E_1}+4
    \,i{\omega}^{3}{ \hat E_2}-8\,{{ \|B^{(0)}\|}}^{2}\omega\,{ \hat E_1}
    \\
    \displaystyle
    & + & 8\,i{{ \|B^{(0)}\| }}^{2}\omega\,{ \hat E_2} ) , 
    \\  
    \displaystyle
    \frac{\hat f_{2}} {f^{(0)}(\|{\bf v}\|^2) \hat D} & = & 
    ik_1{{ v_{\perp}}}^{2} ( -4\,i{k_1}^{2}{{ v_{\perp}}}^{2}{ B^{(0)}}\,{ \hat E_2}+{
      k_1}^{2}{{ v_{\perp}}}^{2}\omega\,{ \hat E_1}-3\,i{k_1}^{2}{{ v_{\perp}}}^{2}\omega\,
    { \hat E_2}-12\,{\omega}^{2}{ B^{(0)}}\,{ \hat E_1}
    \\
    \displaystyle & + & 12\,i{\omega}^{2}{ B^{(0)}}
    \,{ \hat E_2} - 4\,{\omega}^{3}{ \hat E_1}+4\,i{\omega}^{3}{ \hat E_2}-8\,{{ 
        \|B^{(0)}\|}}^{2}\omega\,{ \hat E_1}+8\,i{{ \|B^{(0)}\|}}^{2}\omega\,{ \hat E_2}
    ) ,
  \end{array}
\end{displaymath}
where
\begin{displaymath}
  \hat D = \omega\, ( 64\,{{ \|B^{(0)}\|}}^{4}-16\,{k_1}^{2}{{
  v_{\perp}}}^{2}{\omega}^{ 2}+16\,{\omega}^{4}+16\,{k_1}^{2}{{
  v_{\perp}}}^{2}{{ \|B^{(0)}\|}}^{2}+3\,{k_1}^{4} {{
  v_{\perp}}}^{4}-80\,{{ \|B^{(0)}\|}}^{2}{\omega}^{2} ) .
\end{displaymath}
$k_1 v_{\perp}$ being small with respect to $B^{(0)}$ and $\omega$,
powers of $k_1 v_{\perp}$ can be neglected compared to these terms. The
solution can be written
\begin{displaymath}
  \begin{array}{rcl}
    \displaystyle
    \frac{\hat f_{-2}} {f^{(0)}(\|{\bf v}\|^2) \hat D} & = & 
    i {{v_{\perp}}}^{2} k_1 ( -4\,{\omega}^{3}{ \hat E_1}-4\,i{\omega}^{3}{ \hat E_2}+12\,i{
      \omega}^{2}{ B^{(0)}}\,{ \hat E_2} 
    \\
    \displaystyle 
     & + & 12\,{\omega}^{2}{ B^{(0)}}\,{ \hat E_1}-8\,{{
        \|B^{(0)}\|}}^{2}\omega\,{ \hat E_1}
    - 8\,i{{ \|B^{(0)}\|}}^{2}\omega\,{\hat E_2}), 
    \\
    \displaystyle
    \frac{\hat f_{-1}} {f^{(0)}(\|{\bf v}\|^2) \hat D} & = & 
    2\,i{ v_{\perp}}\, (4\,i
    { B^{(0)}}\,{\omega}^{3}{ \hat E_2}-16\,{{ \|B^{(0)}\|}}^{3}\omega\,{ \hat E_1} - 16\, i{{ \|B^{(0)}\|}}^{3}\omega\,{ \hat E_2}
    - 4\,{ \hat E_1}\,{\omega}^{4}
    \\ 
    \displaystyle
    & + & 16\,{ \hat E_1}\,{{ \|B^{(0)}\|}}^{2}{\omega}^{2} + 16\,i{ \hat E_2}\,{{ \|B^{(0)}\|}}^{2}{\omega}^{2}+4\,
    { B^{(0)}}\,{\omega}^{3}{ \hat E_1}-4\,i{ \hat E_2}\,{\omega}^{4} 
    ),
    \\
    \displaystyle
    \frac{\hat f_{0}} {f^{(0)}(\|{\bf v}\|^2) \hat D} & = & 
    2\,i{{ v_{\perp}}}^{2}k_1 ( 16\,{{ \|B^{(0)}\|}}^{2}\omega\,{ \hat E_1}+{k_1}^{2}{
      { v_{\perp}}}^{2}\omega\,{ \hat E_1}-4\,{\omega}^{3}{ \hat E_1}+4\,i{\omega}^{2}
    { B^{(0)}}\,{ \hat E_2}
    \\
    \displaystyle
    & - & 16\,i{{ \|B^{(0)}\|}}^{3}{ \hat E_2}
    ) ,  
    \\ 
    \displaystyle
    \frac{\hat f_{1}} {f^{(0)}(\|{\bf v}\|^2) \hat D} & = & 
    2\,i { v_{\perp}}\,( \omega-2\,{ B^{(0)}} ) ( -12\,{\omega}^{2}{ B^{(0)}}\,{
      \hat E_1}
    + 12\,i{\omega}^{2}{ B^{(0)}}\,{ \hat E_2}-4\,{\omega}^{3}{ \hat E_1}+4
    \,i{\omega}^{3}{ \hat E_2}
    \\
    \displaystyle
    & - & 
    8\,{{ \|B^{(0)}\|}}^{2}\omega\,{ \hat E_1}+8\,i{{ \|B^{(0)}\|
      }}^{2}\omega\,{ \hat E_2} ) , 
    \\  
    \displaystyle
    \frac{\hat f_{2}} {f^{(0)}(\|{\bf v}\|^2) \hat D} & = & 
    ik_1{{ v_{\perp}}}^{2} (-12\,{\omega}^{2}{ B^{(0))}}\,{ \hat E_1}+12\,i{\omega}^{2}{ B^{(0)}}
    \,{ \hat E_2}
    - 4\,{\omega}^{3}{ \hat E_1}+4\,i{\omega}^{3}{ \hat E_2} 
    \\
    \displaystyle
    & - & 8\,{{ \|B^{(0)}\|}}^{2}\omega\,{ \hat E_1}+8\,i{{ \|B^{(0)}\|}}^{2}\omega\,{ \hat E_2}
    ) ,
  \end{array}
\end{displaymath}
where
\begin{displaymath}
  \hat D = \omega\, ( 64\,{{ \|B^{(0)}\|}}^{4}+16\,{\omega}^{4}-80\,{{
  \|B^{(0)}\|}}^{2}{\omega}^{2} ) .
\end{displaymath}
We choose to initialise the perturbation from the amplitude of the magnetic field:
\begin{displaymath}
  \hat B_3 = A \text{ where } A \in [0,1]. 
\end{displaymath}
Then from the system \eqref{sys_dispersif} and the dispersion relation
\eqref{testtest1}, we deduce the values of $\hat E_1$, $\hat E_2$ and
thus reconstruct the $\hat f_i$,
\begin{displaymath}
  \hat E_1 = {\frac {-i{\hat B_3}\, \left( {\omega}^{4}{\beta}^{2}-{\omega}^{2}{k_1}^{2
        }-{\omega}^{2}{\beta}^{2}-{{ \|B^{(0)}\|}}^{2}{\omega}^{2}{\beta}^{2}+{{ 
            \|B^{(0)}\|}}^{2}{k_1}^{2} \right) }{k_1{\beta}^{2}{ B^{(0)}}}} ,
  \, \, \, 
  \hat E_2 = \frac{\omega \hat B_3}{k_1} .
\end{displaymath}

\end{document}